\documentclass{article}
\usepackage{pdfpages}

\def\bkR{{\rm I\kern-.17em R}}
\def\bkC{{\rm ^{_|}\kern-.47em C}}
\begin{document}
\title{ }
\author{}
\date{}
\begin{center}
{\large\bf Comments on the estimate for Pareto Distribution }\\
\vspace{1cm} U. J. Dixit {\footnote {E-mail:
ulhasdixit@yahoo.co.in} } \\
M. Jabbari Nooghabi {\footnote {E-mail: jabbarinm@yahoo.com, jabbarinm@um.ac.ir} } \\
 {\it Department of Statistics, University of Mumbai, Mumbai-India}\\
\vspace{.1cm}
\end{center} \vspace{.5cm}

\begin{abstract}
Dixit and Jabbari Nooghabi (2010) had derived the MLE and UMVUE of the probability density function (pdf) and cumulative distributive function (cdf). Further, it had been shown that MLE is more efficient than UMVUE. He, Zhou and Zhang (2014) have also derived the same and made a remark that the work of Dixit and Jabbari Nooghabi (2010) is not correct. We have made a comments with detail algebra that our results are correct. Further, we have also given the R code.
\end{abstract}


\section{Introduction}
~~~~The Pareto distribution has been used in connection with studies of income, property values, insurance risk, migration, size of cities and firms, word frequencies, business mortality, service time in queuing systems, etc.\\
The objective of this paper is to discuss efficient estimation of pdf and CDF of Pareto distribution which has been one of the most distinguished candidates for the honor of explaining the distribution of incomes, assets, etc.\\
We assume that random variable $X$ has Pareto distribution with parameter $\alpha$ and $k$ (known) and its probability density function (pdf) is as,
$$ f_X(x)=\frac{\alpha k^\alpha}{x^{\alpha+1}},~~~~~~0<k\leq x,~~ \alpha>0$$
and distribution function
$$F_X(x)=1-\left(\frac{k}{x}\right)^\alpha,~~~~~k\leq x$$
In economics, where this distribution is used as an income distribution, $k$ is some minimum income with a known value. Asrabadi (1990) derived the uniformly minimum variance unbiased estimator (UMVUE) of the probability density function (pdf), the distribution function (cdf) and the $r^{th}$ moment.\\
In this paper, we will give the detail algebra of Dixit and Jabbari Nooghabi's (2010) paper. Also, We have made a comments that our results are correct. Further, we have also given the R code.\\
Dixit and Jabbari Nooghabi (2010) had derived the MLE and UMVUE of the probability density function (pdf) and cumulative distributive function (cdf). Further, it had been shown that MLE is more efficient than UMVUE.
He, Zhou and Zhang (2014) have also derived the same and made a remark that the work of Dixit and Jabbari Nooghabi (2010) is not correct.\\
We like to make some comments as follows.\\
1. We have verified our results and they are correct. We have given the detail algebra and R program. See the attachment.\\
2. Examples given by He et al. (2014) are not correct. One should note MSE of $\hat{f}(x)$ or $\hat{F}(x)$ is a function of parameters. By this one cannot prove anything. Only, one can calculate $\hat{f}(x)$ or $\hat{F}(x)$.\\
3. According to definition of the modified Bessel function in Olver, Lozier, Boisvert, et al. (2010) in the Theorem 1. the notation $K_{(n-r)}(2\sqrt{nr\alpha z})$ should be $K_{(r-n)}(2\sqrt{nr\alpha z})$. Also, the notation $K_{(n)}(2\sqrt{nr\alpha z})$ should be $K_{(-n)}(2\sqrt{nr\alpha z})$.\\
4. After the Theorem 2. the Kummer confluent hypergeometric function is wrong and the correct version is
$$U(a,b,c)=\frac{1}{\Gamma(a)}=\int_0^{\infty}t^{a-1}(1+t)^{b-a-1}e^{-ct}dt.$$

\section{Main Result}
In this section, we give the detail algebra of the paper Dixit and Jabbari Nooghabi (2010).\\
The details of finding result of the second chapter of that paper are as follows.
\begin{eqnarray}
X_{1},\ldots, X_{n} \stackrel{iid}{\sim} f(x)=\frac{\alpha k^\alpha}{x^{\alpha+1}}, ~~~ 0<k\leq x,~\alpha >0
\end{eqnarray}

\begin{eqnarray}
L(x_1,...,x_n,\alpha,k)=\frac{\alpha^nk^{n\alpha}}{\prod_{i=1}^nx_i^{\alpha +1}}\prod_{i=1}^n \textbf {I}(x_i-k)\nonumber,
\end{eqnarray}
where ${\mathbf I}$ is the indicator function defined as
\begin{eqnarray}
{\textbf I}(y) = \left\{
\begin{array}{ll}
1 & \mbox{ $ y>0$},\\
0 &\mbox{ $ otherwise. $}
\end{array}
\right.\nonumber
\end{eqnarray}\begin{eqnarray}
\Rightarrow \ln L(\underline{x},\alpha ,k)=n\ln(\alpha) +n\alpha \ln(k)-(\alpha +1)\ln(\sum_{i=1}^nx_i)\nonumber\\
\Rightarrow \frac{\partial \ln L(\underline{x},\alpha)}{\partial \alpha}=\frac{n}{\alpha}+n\ln(k)-\ln\sum_{i=1}^nx_i=0\nonumber\\
\Rightarrow \tilde{\alpha}= MLE(\alpha)=\frac{n}{\sum_{i=1}^n\ln(\frac{x_i}{k})}, ~~~\frac{\partial^2 \ln L}{\partial\alpha^2}=\frac{-n}{\alpha^2}<0 \nonumber\\
\Rightarrow MLE~of~f(x)=f(x,\tilde{\alpha})=\frac{\tilde{\alpha}k^{\tilde{\alpha}}}{x^{\tilde{\alpha}+1}}\nonumber\\
\Rightarrow \tilde{f}(x)=\frac{\tilde{\alpha}k^{\tilde{\alpha}}}{x^{\tilde{\alpha}+1}},~~~\tilde{\alpha}>0,~0<k\leq x,\nonumber\\
\Rightarrow \tilde{F}(x)=1-\left(\frac{k}x\right)^{\tilde{\alpha}},~~~\tilde{\alpha}>0,~0<k\leq x.
\end{eqnarray}

Put:
\begin{eqnarray}
y=\ln\left(\frac{x}{k}\right)\Rightarrow dy=\frac{1}{x}dx\nonumber\\
\stackrel{x=ke^y}{\Longrightarrow} f_Y(y)=ke^y\frac{\alpha k^\alpha}{(ke^y)^{\alpha+1}}=\alpha e^{-\alpha y},~~y>0,~x\geq k\nonumber\\
\Rightarrow f_Y(y)=\alpha e^{-\alpha y},~~y>0  ~ or  ~ Y{\sim}\Gamma(1,\frac{1}{\alpha})\nonumber\\
\Rightarrow S=\sum_{i=1}^nY_i{\sim}\Gamma(n,\frac{1}{\alpha}),\nonumber\\
g_S(s)=\frac{\alpha^ns^{n-1}}{\Gamma(n)}\exp(-\alpha s),~s>0.
\end{eqnarray}
Let $w=\tilde{\alpha}$
\begin{eqnarray}
(1)\Rightarrow w=\frac{n}{s}\Rightarrow s=\frac{n}{w}\Rightarrow \frac{ds}{dw}=\frac{-n}{w^2}\nonumber\\
\Rightarrow Jacobian=\vert \frac{-n}{w^2} \vert=\frac{n}{w^2}\nonumber
{\Longrightarrow} g(w)=\frac{n}{w^2}\frac{1}{\Gamma(n)(\frac{1}{\alpha})^n}(\frac{n}{w})^{n-1}e^{-n\alpha/w},~~w>0,~\alpha>0\nonumber\\
\Rightarrow g(w)=\frac{(\alpha n)^n}{\Gamma(n)(w^{n+1})}\exp\left\{-\frac{\alpha n}w\right\}, ~w>0.
\end{eqnarray}

\begin{eqnarray}
E(\tilde{\alpha})=E(W)=\int_0^\infty wg(w)dw=\int_0^\infty w\frac{(\alpha n)^n}{\Gamma(n)(w^{n+1})}\exp\left\{-\frac{\alpha n}w\right\}dw=\frac{(\alpha n)^n}{\Gamma(n)}\int_0^\infty w^{-n}e^{-\frac{\alpha n}w}dw.\nonumber
\end{eqnarray}
Put $z=\frac{1}w\Rightarrow dz=-\frac{1}{w^2}dw$, so
\begin{eqnarray}
E(W)&=&\frac{(\alpha n)^n}{\Gamma(n)}\int_0^\infty z^{n-2}e^{-\alpha nz}dz=\frac{(\alpha n)^n}{\Gamma(n)}\Gamma(n-1)(\frac{1}{\alpha n})^{n-1}\int_0^\infty \frac{z^{n-2}e^{-\alpha nz}}{\Gamma(n-1)(\frac{1}{\alpha n})^{n-1}}dz.\nonumber\\
\Rightarrow E(\tilde{\alpha})&=&\frac{\alpha n}{n-1}.\nonumber
\end{eqnarray}
\begin{eqnarray}
E(\tilde{\alpha}^2)=E(W^2)=\int_0^\infty w^2g(w)dw=\int_0^\infty w^2\frac{(\alpha n)^n}{\Gamma(n)(w^{n+1})}\exp\left\{-\frac{\alpha n}w\right\}dw=\frac{(\alpha n)^n}{\Gamma(n)}\int_0^\infty w^{-n+1}e^{-\frac{\alpha n}w}dw.\nonumber
\end{eqnarray}
Same as the previous
\begin{eqnarray}
E(W^2)&=&\frac{(\alpha n)^n}{\Gamma(n)}\int_0^\infty z^{n-3}e^{-\alpha nz}dz=\frac{(\alpha n)^n}{\Gamma(n)}\Gamma(n-2)(\frac{1}{\alpha n})^{n-2}\int_0^\infty \frac{z^{n-3}e^{-\alpha nz}}{\Gamma(n-2)(\frac{1}{\alpha n})^{n-2}}dz.\nonumber\\
\Rightarrow E(\tilde{\alpha}^2)&=&\frac{(\alpha n)^2}{(n-1)(n-2)}.\nonumber
\end{eqnarray}
Therefore
\begin{eqnarray}
MSE(W)&=&V(W)+(E(W)-\alpha)^2=E(W^2)-2\alpha E(W)+\alpha^2\nonumber\\&=&\frac{(\alpha n)^2}{(n-1)(n-2)}-2\alpha\frac{\alpha n}{n-1}+\alpha^2.\nonumber \\ \Rightarrow MSE(\tilde{\alpha})&=&MSE(W)=\frac{\alpha^2(n^2+n-2)}{(n-1)^2(n-2)}.\nonumber
\end{eqnarray}
\textbf{Proof of the Theorem 1.:}\\
(A)
\begin{eqnarray}
E(\tilde{f}(x))&=&\int\tilde{f}(x)g(w)dw
=\int _0^\infty \frac{w k^w}{x^{w+1}}\frac{(\alpha n)^n}{\Gamma(n)w^{n+1}}e^{-\alpha n/w}dw\nonumber\\
&=&\frac{(\alpha n)^n}{\Gamma(n)x}\int _0^\infty\left(\frac{k}{x}\right)^w \frac{e^{-\alpha n/w}}{w ^n}dw.\nonumber
\end{eqnarray}
Put $(\frac{k}{x})^w = e^{w \ln(\frac{k}{x})}$, then
\begin{eqnarray}
E(\tilde{f}(x))&=&\frac{(\alpha n)^n}{\Gamma(n)x}\int _0^\infty\frac{e^{w ln(\frac{k}{x})}e^{-\alpha n/w}}{w^n}dw.\nonumber
\end{eqnarray}
We know $e^{w \ln(\frac{k}{x})}=\sum_{j=0}^\infty \frac{w^j(\ln\frac{k}{x})^j}{j!}$
\begin{eqnarray}
E(\tilde{f}(x))&=&\frac{(\alpha n)^n}{\Gamma(n)x}\int _0^\infty\sum_{j=0}^\infty \frac{w^j(\ln(\frac{k}{x}))^j}{j! w^n}e^{-\alpha n/w}dw\nonumber\\
&=&\frac{(\alpha n)^n}{\Gamma(n)x}\sum_{j=0}^\infty\frac{(\ln(\frac{k}{x}))^j}{j!}\int _0^\infty\frac{e^{-\alpha n/w}}{w^{n-j}}dw.\nonumber
\end{eqnarray}
Put $\frac{1}{w}=z$ then $\frac{-1}{w^2}dw=dz\Rightarrow dw=\frac{-dz}{z^2}$.
\begin{eqnarray}
E(\tilde{f}(x))&=&\frac{(\alpha n)^n}{\Gamma(n)x}\sum_{j=0}^\infty\frac{(\ln(\frac{k}{x}))^j}{j!}\int _0^\infty z^{n-j-2}e^{-\alpha nz}dz\nonumber\\
&=&\frac{(\alpha n)^n}{\Gamma(n)x}\sum_{j=0}^\infty\frac{(\ln(\frac{k}{x}))^j}{j!}\frac{\Gamma(n-j-1)}{(\frac{1}{n\alpha})^{-n+j+1}}\int _0^\infty\frac{z^{n-j-2}e^{-\alpha nz}}{\Gamma(n-j-1)(\frac{1}{n\alpha})^{n-j-1}}dz\nonumber\\
&=&\frac{1}{\Gamma(n)x}\sum_{j=0}^\infty\frac{(n\alpha)^{j+1}}{j!}\Gamma(n-j-1)\left(\ln\left(\frac{k}{x}\right)\right)^j\nonumber\\
\Rightarrow E(\tilde{f}(x))&=&\frac{1}{\Gamma(n)x}\sum_{j=0}^{n-2}\frac{(n\alpha)^{j+1}}{j!}\Gamma(n-j-1)\left(\ln\left(\frac{k}{x}\right)\right)^j.
\end{eqnarray}
(B)

\begin{eqnarray}
E(\tilde{F}(x))&=&\int\tilde{F}(x)g(w)dw=\int_0^\infty\left[1-(\frac{k}{x})^w\right]\frac{(\alpha n)^n}{\Gamma(n)w^{n+1}}e^{-\alpha n/w}dw\nonumber\\
&=&\int_0^\infty\frac{(\alpha n)^n}{\Gamma(n)w^{n+1}}e^{-\alpha n/w}dw-\frac{(\alpha n)^n}{\Gamma(n)}\int_0^\infty(\frac{k}{x})^w \frac{1}{w^{n+1}}e^{-\alpha n/w}dw\nonumber\\
&=&1-\frac{(\alpha n)^n}{\Gamma(n)}\int_0^\infty\frac{e^{w \ln\left(\frac{k}{x}\right)}e^{-\alpha n/w}}{w^{n+1}}dw.\nonumber
\end{eqnarray}
We know that $\left(\frac{k}{x}\right)^w =e^{w \ln\left(\frac{k}{x}\right)}=\sum_{j=0}^\infty\frac{w^j(\ln\left(\frac{k}{x}\right))^j}{j!}$, then
\begin{eqnarray}
E(\tilde{F}(x))&=&1-\frac{(\alpha n)^n}{\Gamma(n)}\int_0^\infty\sum_{j=0}^\infty\frac{w^j(\ln(\frac{k}{x}))^j}{j!}w^{-n-1}e^{-\alpha n/w}dw\nonumber\\
&=&1-\frac{(\alpha n)^n}{\Gamma(n)}\sum_{j=0}^\infty\frac{(\ln(\frac{k}{x}))^j}{j!}\int_0^\infty w^{-n-1+j}e^{-\alpha n/w}dw\nonumber\\
&=&1-\frac{(\alpha n)^n}{\Gamma(n)}\sum_{j=0}^\infty\frac{(\ln(\frac{k}{x}))^j}{j!}\int_0^\infty\frac{e^{-\alpha n/w}}{w^{n+1-j}}dw\nonumber.
\end{eqnarray}
Put $\frac{1}{w}=z\Rightarrow\frac{-1}{w^2}dw=dz \Rightarrow dw=\frac{-1}{z^2}dz.$
Then
\begin{eqnarray}
E(\tilde{F}(x))&=&1-\frac{(\alpha n)^n}{\Gamma(n)}\sum_{j=0}^\infty\frac{(\ln(\frac{k}{x}))^j}{j!}\Gamma(n-j)\left(\frac{1}{\alpha n}\right)^{n-j}\int_0^\infty\frac{z^{n-j-1}e^{-\alpha nz}}{\Gamma(n-j)(\frac{1}{\alpha n})^{n-j}}dz\nonumber\\
\Rightarrow E(\tilde{F}(x))&=&1-\frac{1}{\Gamma(n)}\sum_{j=0}^{n-1}\frac{(\alpha n)^j}{j!}\Gamma(n-j)\left(\ln\left(\frac{k}{x}\right)\right)^j,~~~x\geq k.
\end{eqnarray}
\textbf{Proof of the Theorem 2.:}\\
(A)
At first, we should find $E(\tilde{f}(x))^2$. So
\begin{eqnarray}
E(\tilde{f}(x))^2&=&\int_0^\infty(\tilde{f}(x))^2g(w)dw=\int_0^\infty\frac{w^2k^{2w}}{x^{2w +2}}\frac{(\alpha n)^n}{\Gamma(n)w^{n+1}}e^{-\alpha n/w}dw\nonumber\\
&=&\frac{(\alpha n)^n}{\Gamma(n)x^2}\int_0^\infty \left(\frac{k}{x}\right)^{2w}\frac{e^{-\alpha n/w}}{w^{n-1}}dw.\nonumber
\end{eqnarray}
Similarly to the pervious Theorem, we have
\begin{eqnarray}
E(\tilde{f}(x))^2&=&\frac{(\alpha n)^n}{\Gamma(n)x^2}\int _0^\infty\frac{e^{2w \ln\left(\frac{k}{x}\right)}e^{-\alpha n/w}}{w^{n-1}}dw\nonumber\\
&=&\frac{(\alpha n)^n}{\Gamma(n)x^2}\int _0^\infty\sum_{j=0}^\infty \frac{2^jw^j(\ln\left(\frac{k}{x}\right))^j}{j! w^{n-1}}e^{-\alpha n/w}dw\nonumber\\
&=&\frac{(\alpha n)^n}{\Gamma(n)x^2}\sum_{j=0}^\infty\frac{2^j(\ln\left(\frac{k}{x}\right))^j}{j!}\int _0^\infty\frac{e^{-\alpha n/w}}{w^{n-j-1}}dw\nonumber\\
&\stackrel{\frac{1}{w}=z}{=}&\frac{(\alpha n)^n}{\Gamma(n)x^2}\sum_{j=0}^\infty\frac{2^j(\ln\left(\frac{k}{x}\right))^j}{j!}\int _0^\infty z^{n-j-3}e^{-\alpha nz}dz\nonumber\\
&=&\frac{(\alpha n)^n}{\Gamma(n)x^2}\sum_{j=0}^\infty\frac{2^j(\ln\left(\frac{k}{x}\right))^j}{j!}\Gamma(n-j-2)\left(\frac{1}{\alpha n}\right)^{n-j-2}\int _0^\infty\frac{z^{n-j-3}e^{-\alpha nz}}{\Gamma(n-j-2)(\frac{1}{\alpha n})^{n-j-2}}dz\nonumber\\
\Rightarrow E(\tilde{f}(x))^2&=&\frac{1}{\Gamma(n)x^2}\sum_{j=0}^{n-2}\frac{2^j(\ln\left(\frac{k}{x}\right))^j}{j!}\Gamma(n-j-2)(\alpha n)^{j+2}.\nonumber
\end{eqnarray}
We know that $V(\tilde{f}(x))=E(\tilde{f}(x))^2-E^2(\tilde{f}(x))$. Then
\begin{eqnarray}
V(\tilde{f}(x))&=&\frac{1}{\Gamma(n)x^2}\sum_{j=0}^{n-2}\frac{2^j(\ln\left(\frac{k}{x}\right))^j}{j!}\Gamma(n-j-2)(\alpha n)^{j+2}\nonumber\\
&-&\left[\frac{1}{\Gamma(n)x}\sum_{j=0}^{n-1}\frac{(\alpha n)^{j+1}}{j!}\Gamma(n-j-1)(\ln\left(\frac{k}{x}\right))^j\right]^2\nonumber.
\end{eqnarray}
Therefore
\begin{eqnarray}
MSE(\tilde{f}(x))&=&V(\tilde{f}(x))+(E(\tilde{f}(x))-f(x))^2\nonumber\\
&=&V(\tilde{f}(x))+E^2(\tilde{f}(x))-2E(\tilde{f}(x))f(x)+f^2(x)\nonumber\\
&=&E(\tilde{f}(x))^2-E^2(\tilde{f}(x))+E^2(\tilde{f}(x))-2E(\tilde{f}(x))f(x)+f^2(x)\nonumber\\
&=&E(\tilde{f}(x))^2-2f(x)E(\tilde{f}(x))+f^2(x)\nonumber\\
\Rightarrow MSE(\tilde{f}(x))&=&\frac{1}{\Gamma(n)x^2}\sum_{j=0}^{n-2}\frac{2^j(\ln\left(\frac{k}{x}\right))^j}{j!}\Gamma(n-j-2)(\alpha n)^{j+2}\nonumber\\&-&2\frac{\alpha k^\alpha}{x^{\alpha+1}}\frac{1}{\Gamma(n)x}\sum_{j=0}^{n-1}\frac{(n\alpha)^{j+1}}{j!}\Gamma(n-j-1)
(\ln\left(\frac{k}{x}\right))^j+(\frac{\alpha k^\alpha}{x^{\alpha+1}})^2.
\end{eqnarray}
(B) Seme as the case (A)
\begin{eqnarray}
E(\tilde{F}(x))^2&=&\int_0^\infty(\tilde{F}(x))^2g(w)dw=\int_0^\infty\left[1-\left(\frac{k^w}{x^w}\right)\right]^2\frac{(\alpha n)^n}{\Gamma(n)w^{n+1}}e^{-\alpha n/w}dw\nonumber\\
&=&\int_0^\infty\frac{(\alpha n)^n}{\Gamma(n)w^{n+1}}e^{-\alpha n/w}dw-2\frac{(\alpha n)^n}{\Gamma(n)}\int_0^\infty\frac{(\frac{k}{x})^w e^{-\alpha n/w}}{w^{n+1}}dw+\frac{(\alpha n)^n}{\Gamma(n)}\int_0^\infty\frac{(\frac{k}{x})^{2w }e^{-\alpha n/w}}{w^{n+1}}dw \nonumber\\
&=&1-2\frac{(\alpha n)^n}{\Gamma(n)}\int_0^\infty\frac{e^{w \ln\left(\frac{k}{x}\right)} e^{-\alpha n/w}}{w^{n+1}}dw+\frac{(\alpha n)^n}{\Gamma(n)}\int_0^\infty\frac{e^{2w \ln\frac{k}{x}}e^{-\alpha n/w}}{w^{n+1}}dw.\nonumber
\end{eqnarray}
Let $e^{w \ln\left(\frac{k}{x}\right)}=\sum_{j=0}^\infty\frac{w^j(\ln\left(\frac{k}{x}\right))^j}{j!}, ~~e^{2w
\ln\left(\frac{k}{x}\right)}=\sum_{j=0}^\infty\frac{2^jw^j(\ln\left(\frac{k}{x}\right))^j}{j!}$, then
\begin{eqnarray}
E(\tilde{F}(x))^2&=&1-2\frac{(\alpha n)^n}{\Gamma(n)}\int_0^\infty\sum_{j=0}^\infty\frac{w^j(\ln\left(\frac{k}{x}\right))^j}{j!}\frac{e^{-\alpha n/w}}{w^{n+1}}dw \nonumber\\
&+&\frac{(\alpha n)^n}{\Gamma(n)}\int_0^\infty\sum_{j=0}^\infty\frac{2^jw^j(\ln\left(\frac{k}{x}\right))^j}{j!}\frac{e^{-\alpha n/w}}{w^{n+1}}dw\nonumber\\
&=&1-2\frac{(\alpha n)^n}{\Gamma(n)}\sum_{j=0}^\infty\frac{(\ln\left(\frac{k}{x}\right))^j}{j!}\int_0^\infty\frac{e^{-\alpha n/w}}{w^{n+1-j}}dw\nonumber\\
&+&\frac{(\alpha n)^n}{\Gamma(n)}\sum_{j=0}^\infty\frac{2^j(\ln\left(\frac{k}{x}\right))^j}{j!}\int_0^\infty\frac{e^{-\alpha n/w}}{w^{n+1-j}}dw\nonumber\\
&\stackrel{\frac{1}{w}=z}{=}&
1-2\frac{(\alpha n)^n}{\Gamma(n)}\sum_{j=0}^\infty\frac{(\ln\left(\frac{k}{x}\right))^j}{j!}\Gamma(n-j)\left(\frac{1}{\alpha n}\right)^{n-j}\int_0^\infty\frac{z^{n-1-j}e^{-\alpha nz}}{\Gamma(n-j)\left(\frac{1}{\alpha n}\right)^{n-j}}dz\nonumber\\
&+&\frac{(\alpha n)^n}{\Gamma(n)}\sum_{j=0}^\infty\frac{2^j(\ln\left(\frac{k}{x}\right))^j}{j!}\Gamma(n-j)\left(\frac{1}{\alpha n}\right)^{n-j}\int_0^\infty\frac{z^{n-1-j}e^{-\alpha nz}}{\Gamma(n-j)\left(\frac{1}{\alpha n}\right)^{n-j}}dz\nonumber\\
\Rightarrow E(\tilde{F}(x))^2 &=&1-\frac{2}{\Gamma(n)}\sum_{j=0}^{n-1}\frac{\Gamma(n-j)(\alpha n)^j}{j!}(\ln\left(\frac{k}{x}\right))^j\nonumber\\&+&\frac{1}{\Gamma(n)}\sum_{j=0}^{n-1}\frac{\Gamma(n-j){2^j}(\alpha n)^j}{j!}(\ln\left(\frac{k}{x}\right))^j, ~~~x\geq k.\nonumber
\end{eqnarray}
We have $V(\tilde{F}(x))=E(\tilde{F}(x))^2-E^2(\tilde{F}(x))$, so
\begin{eqnarray}
V(\tilde{F}(x))&=&1-\frac{2}{\Gamma(n)}\sum_{j=0}^{n}\frac{\Gamma(n-j)(\alpha n)^j}{j!}(\ln\left(\frac{k}{x}\right))^j+\frac{1}{\Gamma(n)}\sum_{j=0}^{n}\frac{\Gamma(n-j){2^j}(\alpha n)^j}{j!}(\ln\left(\frac{k}{x}\right))^j\nonumber\\
&-&\left[1-\frac{1}{\Gamma(n)}\sum_{j=0}^{n}\frac{\Gamma(n-j)(\alpha n)^j}{j!}(\ln\left(\frac{k}{x}\right))^j\right]^2\nonumber.
\end{eqnarray}
Further $MSE(\tilde{F}(x))=V(\tilde{F}(x))+(E(\tilde{F}(x))-F(x))^2=E(\tilde{F}(x))^2-2F(x)E(\tilde{F}(x))+F^2(x)$, then
\begin{eqnarray}
MSE(\tilde{F}(x))&=&1-\frac{2}{\Gamma(n)}\sum_{j=0}^{n}\frac{\Gamma(n-j)(\alpha n)^j}{j!}(\ln\left(\frac{k}{x}\right))^j+\frac{1}{\Gamma(n)}\sum_{j=0}^{n}\frac{\Gamma(n-j){2^j}(\alpha n)^j}{j!}(\ln\left(\frac{k}{x}\right))^j\nonumber\\
&-&2\left[1-\left(\frac{k}{x}\right)^\alpha\right]\left[1-\frac{1}{\Gamma(n)}\sum_{j=0}^{n}\frac{\Gamma(n-j)(\alpha n)^j}{j!}(\ln\left(\frac{k}{x}\right))^j\right]+\left[1-\left(\frac{k}{x}\right)^\alpha\right]^2.\nonumber
\end{eqnarray}
Therefore
\begin{eqnarray}
MSE(\tilde{F}(x))&=&2+\frac{1}{\Gamma(n)}\sum_{j=0}^{n}\frac{\Gamma(n-j){2^j}(\alpha n)^j}{j!}(\ln\left(\frac{k}{x}\right))^j\nonumber\\
&-&2\left(\frac{k}{x}\right)^\alpha\frac{1}{\Gamma(n)}\sum_{j=0}^{n}\frac{\Gamma(n-j)(\alpha n)^j}{j!}(\ln\left(\frac{k}{x}\right))^j+\left(\frac{k}{x}\right)^{2\alpha}.
\end{eqnarray}
From Asrabadi (1990), we have
\begin{eqnarray}
\hat{\alpha}=\frac{n-1}{\ln(t)-n\ln(k)},~~t\geq k^n,
\end{eqnarray}
the UMVUE of $f(x)$ and $F(x)$ is
\begin{eqnarray}
\hat{f}(x)=\frac{(n-1)[\ln(t)-\ln(x)-(n-1)\ln(k)]^{n-2}}{x[\ln(t)-n\ln(k)]^{n-1}},~~k\leq x< tk^{1-n},
\end{eqnarray}
\begin{eqnarray}
\hat{F}(x)=\left\{
\begin{array}{ll}
0& \mbox{$x<k$},\\
1-\frac{[\ln(t)-\ln(x)-(n-1)\ln(k)]^{n-1}}{[\ln(t)-n\ln(k)]^{n-1}}& \mbox{$k\leq x\leq tk^{1-n}$},\\
1& \mbox{$x\geq tk^{1-n}$},
\end{array}
\right.
\end{eqnarray}
respectively. Also,
\begin{eqnarray}
f(x)=\frac{\alpha k^\alpha}{x^{\alpha+1}},~~ x\geq k>0,~\alpha>0\nonumber,
\end{eqnarray}
\begin{eqnarray}
F(x)=1-\left(\frac{k}x\right)^\alpha,~~k\leq x\nonumber.
\end{eqnarray}
\textbf{Proof of the Theorem 3.:}\\
(A)
It is obvious that $E(\hat{f}(x))=f(x)$. So, we should find
\begin{eqnarray}
E(\hat{f}(x))^2= \int \hat{f}^2(x)h^*(x)dt\nonumber,
\end{eqnarray}
where,
\begin{eqnarray}
h^*(x)=\frac{\alpha^n k^{n\alpha}}{(n-1)!}t^{-\alpha-1}[\ln(t)-n\ln(k)]^{n-1},~~t\geq k^n\nonumber.
\end{eqnarray}
Therefore
\begin{eqnarray}
E(\hat{f}(x))^2&=& \int_{xk^{n-1}}^\infty\frac{(n-1)^2[\ln(t)-\ln(x)-(n-1)\ln(k)]^{2n-4}}{x^2[\ln(t)-n\ln(k)]^{2n-2}}\frac{\alpha^n k^{n\alpha}}{(n-1)!}t^{-\alpha-1}[\ln(t)-n\ln(k)]^{n-1}dt\nonumber\\
&=&\frac{(n-1)^2\alpha^n k^{\alpha n}}{x^2(n-1)!}\int_{xk^{n-1}}^\infty\frac{[\ln(t)-\ln(x)-(n-1)\ln(k)]^{2n-4}}{[\ln(t)-n\ln(k)]^{n-1}}t^{-\alpha-1}dt.\nonumber
\end{eqnarray}
Let $z=\ln(t)-n\ln(k)\Rightarrow dz=\frac{1}{t}dt$, then
\begin{eqnarray}
E(\hat{f}(x))^2=\frac{(n-1)\alpha^n}{x^2(n-2)!}\int_{\ln\left(\frac{x}{k}\right)}^\infty \frac{[z-\ln(x)+\ln(k)]^{2n-4}e^{-\alpha z}}{z^{n-1}}dz\nonumber.
\end{eqnarray}
We know that $$[z-\ln(x)+\ln(k)]^{2n-4}=\sum_{j=0}^{2n-4}C(2n-4,j)\left(-\ln\left(\frac{x}k\right)\right)^jz^{2n-4-j},$$
where $C(n,k)=\frac{n!}{k!(n-k)!}$. So
\begin{eqnarray}
E(\hat{f}(x))^2=\frac{(n-1)\alpha^n}{x^2(n-2)!}\sum_{j=0}^{2n-4}C(2n-4,j)\left(-\ln\left(\frac{x}k\right)\right)^j
\int_{\ln\left(\frac{x}{k}\right)}^\infty z^{n-3-j}e^{-\alpha z}dz\nonumber.
\end{eqnarray}
The above integral is the incomplete Gamma function, therefore
\begin{eqnarray}
E(\hat{f}(x))^2&=&\frac{(n-1)\alpha^n}{x^2(n-2)!}\sum_{j=0}^{2n-4}C(2n-4,j)\left(-\ln\left(\frac{x}k\right)\right)^j
\nonumber\\&\times&
\frac{\Gamma(n-2-j)}{\alpha^{n-2-j}}\sum_{i=0}^{n-3-j}\frac{\exp\left(-\alpha\ln\left(\frac{x}k\right)\right)
\left(\alpha\ln\left(\frac{x}k\right)\right)^i}{i!}\nonumber\\
&=&\frac{(n-1)\alpha^2 k^\alpha}{x^{\alpha+2}(n-2)!} \sum_{j=0}^{2n-4}C(2n-4,j)\alpha^j\left[-\ln\left(\frac{x}{k}\right)\right]^j\Gamma(n-2-j)
\sum_{i=0}^{n-3-j}\frac{\alpha^i\left(\ln\left(\frac{x}k\right)\right)^i}{i!}\nonumber.
\end{eqnarray}
We know that the Gamma function is defined on the positive value. So
\begin{eqnarray}
E(\hat{f}(x))^2=\frac{(n-1)\alpha^2 k^\alpha}{x^{\alpha+2}(n-2)!} \sum_{j=0}^{n-3}C(2n-4,j)\alpha^j\left[-\ln\left(\frac{x}{k}\right)\right]^j\Gamma(n-2-j)
\sum_{i=0}^{n-3-j}\frac{\alpha^i\left(\ln\left(\frac{x}k\right)\right)^i}{i!}\nonumber.
\end{eqnarray}
Finally
\begin{eqnarray}
MSE(\hat{f}(x))=V(\hat{f}(x))&=&\frac{(n-1)\alpha^2 k^\alpha}{x^{\alpha+2}\Gamma(n-1)} \sum_{j=0}^{n-3}C(2n-4,j)\alpha^j\Gamma(n-j-2)\left(-\ln\left(\frac{x}{k}\right)\right)^j\nonumber\\
&\times&\sum_{i=0}^{n-3-j}\frac{\alpha^i\left(\ln\left(\frac{x}k\right)\right)^i}{i!}
-\left(\frac{\alpha k^{\alpha}}{x^{\alpha+1}}\right)^2.
\end{eqnarray}
(B)

\begin{eqnarray}
E(\hat{F}(x))^2&=&\int\hat{F}^2(x)h^*(t)dt\nonumber\\
&=&\int_{xk^{n-1}}^\infty\left[1-\frac{[\ln(t)-\ln(x)-(n-1)\ln(k)]^{n-1}}{[\ln(t)-n\ln(k)]^{n-1}}\right]^2h^*(t)dt
+\int_{k^n}^{xk^{n-1}}1^2\times h^*(t)dt\nonumber\\
&=&\int_{xk^{n-1}}^\infty h^*(t)dt-2\int_{xk^{n-1}}^\infty\frac{[\ln(t)-\ln(x)-(n-1)\ln(k)]^{n-1}}{[\ln(t)-n\ln(k)]^{n-1}}h^*(t)dt\nonumber\\
&+&\int_{xk^{n-1}}^\infty\frac{[\ln(t)-\ln(x)-(n-1)\ln(k)]^{2n-2}}{[\ln(t)-n\ln(k)]^{2n-2}}h^*(t)dt
+\int_{k^n}^{xk^{n-1}} h^*(t)dt\nonumber\\
&=&\int_{xk^{n-1}}^\infty h^*(t)dt+\int_{k^n}^{xk^{n-1}} h^*(t)dt-2\int_{xk^{n-1}}^\infty\frac{[\ln(t)-\ln(x)-(n-1)\ln(k)]^{n-1}}{[\ln(t)-n\ln(k)]^{n-1}}\nonumber\\
&\times&\frac{\alpha ^nk^{\alpha n}}{(n-1)!}t^{-\alpha-1}[\ln(t)-n\ln(k)]^{n-1}dt\nonumber\\
&+&\int_{xk^{n-1}}^\infty\frac{[\ln(t)-\ln(x)-(n-1)\ln(k)]^{2n-2}}{[\ln(t)-n\ln(k)]^{2n-2}}\frac{\alpha^nk^{\alpha n}}{(n-1)!}t^{-\alpha-1}[\ln(t)-n\ln(k)]^{n-1}dt\nonumber\\
&=&\int_{k^n}^{\infty}h^*(t)dt-2\frac{\alpha^nk^{\alpha n}}{(n-1)!}\int_{xk^{n-1}}^\infty [\ln(t)-\ln(x)-(n-1)\ln(k)]^{n-1}t^{-\alpha-1}dt\nonumber\\
&+&\frac{\alpha^nk^{\alpha n}}{(n-1)!}\int_{xk^{n-1}}^\infty\frac{[\ln(t)-\ln(x)-(n-1)\ln(k)]^{2n-2}}{[\ln(t)-n\ln(k)]^{n-1}}t^{-\alpha-1}dt.\nonumber
\end{eqnarray}
We know that $\int_{k^n}^{\infty}h^*(t)dt=1$. For second part let $z=\ln(t)-\ln(x)-(n-1)\ln(k)$ and to solve the third integral put $z=\ln(t)-n\ln(k)$. Then
\begin{eqnarray}
E(\hat{F}(x))^2&=&1-2\frac{\alpha^nk^{\alpha n}}{(n-1)!}x^{-\alpha}k^{-\alpha(n-1)}\int_0^\infty z^{n-1}e^{-\alpha z}dz\nonumber\\
&+&\frac{\alpha^nk^{\alpha n}}{(n-1)!}k^{-n\alpha}\int_{\ln\left(\frac{x}k\right)}^\infty \frac{[z-\ln(x)+\ln(k)]^{2n-2}}{z^{n-1}}e^{-\alpha z}dz\nonumber\\
&=&1-2\frac{\alpha^n k^{\alpha}}{x^\alpha(n-1)!}\int_0^\infty z^{n-1}e^{-\alpha z}dz+\frac{\alpha^n }{(n-1)!}\int_{\ln\left(\frac{x}k\right)}^\infty \frac{[z-\ln(x)+\ln(k)]^{2n-2}e^{-\alpha z}}{z^{n-1}}dz\nonumber.
\end{eqnarray}
$\int_0^\infty z^{n-1}e^{-\alpha z}dz=\frac{\Gamma(n)}{\alpha^n}$ and for the last integral, we know $[z-\ln(x)+\ln(k)]^{2n-2}=\sum_{j=0}^{2n-2}C(2n-2,j)z^{2n-2-j}\left[-\ln\left(\frac{x}k\right)\right]^j$. Therefore

\begin{eqnarray}
E(\hat{F}(x))^2&=&1-2\frac{k^{\alpha}}{x^\alpha}+\frac{\alpha^n}{(n-1)!}\sum_{j=0}^{2n-2}C(2n-2,j)
\left[-\ln\left(\frac{x}k\right)\right]^j\int_{\ln\left(\frac{x}k\right)}^\infty z^{n-1-j}e^{-\alpha z}{z^{n-1}}dz\nonumber\\
&=&1-2\frac{k^\alpha}{x^\alpha}+\frac{\alpha^n}{(n-1)!}\sum_{j=0}^{2n-2}C(2n-2,j)
\left[-\ln\left(\frac{x}k\right)\right]^j\frac{\Gamma(n-j)}{\alpha^{n-j}}\int_{\ln\left(\frac{x}k\right)}^\infty \frac{\alpha^{n-j}z^{n-1-j}e^{-\alpha z}}{\Gamma(n-j)}dz\nonumber.
\end{eqnarray}
The last integral is the incomplete Gamma function, then
\begin{eqnarray}
E(\hat{F}(x))^2&=&1-2\frac{k^\alpha}{x^\alpha}+\frac{\alpha^n}{(n-1)!}\sum_{j=0}^{2n-2}C(2n-2,j)
\left[-\ln\left(\frac{x}k\right)\right]^j\nonumber\\&\times&\frac{\Gamma(n-j)}{\alpha^{n-j}}\sum_{i=0}^{n-1-j} \frac{e^{-\alpha\ln\left(\frac{x}k\right)}\left[\alpha\ln\left(\frac{x}k\right)\right]^i}{i!}\nonumber.
\end{eqnarray}
The Gamma function is defined over positive value, So
\begin{eqnarray}
E(\hat{F}(x))^2&=&1-2\frac{k^\alpha}{x^\alpha}+\frac{k^\alpha}{x^\alpha(n-1)!}\sum_{j=0}^{n-1}C(2n-2,j)
\left[-\ln\left(\frac{x}k\right)\right]^j\nonumber\\&\times&\alpha^j\Gamma(n-j)\sum_{i=0}^{n-1-j} \frac{\left[\alpha\ln\left(\frac{x}k\right)\right]^i}{i!}\nonumber.
\end{eqnarray}
Then
\begin{eqnarray}
MSE(\hat{F}(x))&=&V(\hat{F}(x))=E(\hat{F}(x))^2-E^2(\hat{F}(x))=E(\hat{F}(x))^2-F^2(x)\nonumber\\
&=&1-2\frac{k^\alpha}{x^\alpha}+\frac{k^\alpha}{x^\alpha(n-1)!}\sum_{j=0}^{n-1}C(2n-2,j)
\left[-\ln\left(\frac{x}k\right)\right]^j\nonumber\\&\times&\alpha^j\Gamma(n-j)\sum_{i=0}^{n-1-j} \frac{\left[\alpha\ln\left(\frac{x}k\right)\right]^i}{i!}-\left[1-\left(\frac{k}x\right)^\alpha\right]^2.\nonumber
\end{eqnarray}
Finally
\begin{eqnarray}
MSE(\hat{F}(x))&=&\frac{k^\alpha}{\Gamma(n)x^\alpha}\sum_{j=0}^{n-1}C(2n-2,j)
\alpha^j\Gamma(n-j)\left[-\ln\left(\frac{x}k\right)\right]^j\nonumber\\&\times&\sum_{i=0}^{n-1-j} \frac{\alpha^i\left[\ln\left(\frac{x}k\right)\right]^i}{i!}-\left(\frac{k}x\right)^{2\alpha}.
\end{eqnarray}

\subsection{The rth estimate of $\tilde{f}(x)$ and $\tilde{F}(x)$}
To find the the rth estimate of $\tilde{f}(x)$, we have
\begin{eqnarray}
E(\tilde{f}(x))^r&=&\int_0^\infty(\tilde{f}(x))^rg(w)dw=\int_0^\infty\frac{w^rk^{rw}}{x^{r(w+1)}}\frac{(\alpha n)^n}{\Gamma(n)w^{n+1}}e^{-\alpha n/w}dw\nonumber\\
&=&\frac{(\alpha n)^n}{\Gamma(n)x^r}\int_0^\infty \left(\frac{k}{x}\right)^{rw}\frac{e^{-\alpha n/w}}{w^{n-r+1}}dw\nonumber\\
&=&\frac{(\alpha n)^n}{\Gamma(n)x^r}\int _0^\infty\frac{e^{rw \ln\left(\frac{k}{x}\right)}e^{-\alpha n/w}}{w^{n-r+1}}dw\nonumber\\
&=&\frac{(\alpha n)^n}{\Gamma(n)x^r}\int _0^\infty\sum_{j=0}^\infty \frac{r^jw^j(\ln\left(\frac{k}{x}\right))^j}{j! w^{n-r+1}}e^{-\alpha n/w}dw\nonumber\\
&=&\frac{(\alpha n)^n}{\Gamma(n)x^r}\sum_{j=0}^\infty\frac{r^j(\ln\left(\frac{k}{x}\right))^j}{j!}\int _0^\infty\frac{e^{-\alpha n/w}}{w^{n-r-j+1}}dw\nonumber\\
&\stackrel{\frac{1}{w}=z}{=}&\frac{(\alpha n)^n}{\Gamma(n)x^r}\sum_{j=0}^\infty\frac{r^j(\ln\left(\frac{k}{x}\right))^j}{j!}\int _0^\infty z^{n-r-j-1}e^{-\alpha nz}dz\nonumber\\
&=&\frac{(\alpha n)^n}{\Gamma(n)x^r}\sum_{j=0}^\infty\frac{r^j(\ln\left(\frac{k}{x}\right))^j}{j!}\Gamma(n-r-j)\left(\frac{1}{\alpha n}\right)^{n-r-j}\nonumber\\
\Rightarrow E(\tilde{f}(x))^r&=&\frac{1}{\Gamma(n)x^r}\sum_{j=0}^{n-r-1}\frac{r^j(\ln\left(\frac{k}{x}\right))^j}{j!}\Gamma(n-r-j)(\alpha n)^{j+r}.
\end{eqnarray}
Also, the rth estimate of $\tilde{F}(x)$ can be found by calculating the following integral.
\begin{eqnarray}
E(\tilde{F}(x))^r&=&\int_0^\infty(\tilde{F}(x))^rg(w)dw=\int_0^\infty\left[1-\left(\frac{k}{x}\right)^w\right]^r\frac{(\alpha n)^n}{\Gamma(n)w^{n+1}}e^{-\alpha n/w}dw\nonumber\\
&=&\int_0^\infty\sum_{j=0}^{r}C(r,j)\left(-\left(\frac{k}x\right)^w\right)^j\frac{(\alpha n)^n}{\Gamma(n)w^{n+1}}e^{-\alpha n/w}dw\nonumber\\
&=&\frac{(\alpha n)^n}{\Gamma(n)}\sum_{j=0}^{r}C(r,j)(-1)^j\int_0^\infty\frac{(\frac{k}{x})^{jw} e^{-\alpha n/w}}{w^{n+1}}dw\nonumber\\
&=&\frac{(\alpha n)^n}{\Gamma(n)}\sum_{j=0}^{r}C(r,j)(-1)^j\int_0^\infty\frac{e^{jw \ln\left(\frac{k}{x}\right)} e^{-\alpha n/w}}{w^{n+1}}dw\nonumber\\
&=&\frac{(\alpha n)^n}{\Gamma(n)}\sum_{j=0}^{r}C(r,j)(-1)^j\sum_{i=0}^{\infty}\frac{\left(j\ln\left(\frac{k}x\right)\right)^i}{i!}\int_0^\infty \frac{e^{-\alpha n/w}}{w^{n+1-i}}dw\nonumber\\
&\stackrel{\frac{1}{w}=z}{=}&\frac{(\alpha n)^n}{\Gamma(n)}\sum_{j=0}^{r}C(r,j)(-1)^j\sum_{i=0}^{\infty}\frac{\left(j\ln\left(\frac{k}x\right)\right)^i}{i!}\int_0^\infty z^{n-i-1}e^{-\alpha nz}dz\nonumber\\
&=&\frac{(\alpha n)^n}{\Gamma(n)}\sum_{j=0}^{r}C(r,j)(-1)^j\sum_{i=0}^{\infty}\frac{\left(j\ln\left(\frac{k}x\right)\right)^i}{i!} \Gamma(n-i)(\frac{1}{\alpha n})^{n-i}.\nonumber
\end{eqnarray}
Then
\begin{eqnarray}
E(\tilde{F}(x))^r=\frac{1}{\Gamma(n)}\sum_{j=0}^{r}C(r,j)(-1)^j\sum_{i=0}^{n-1}\frac{\left(j\ln\left(\frac{k}x\right)\right)^i}{i!} \Gamma(n-i)(\alpha n)^{i}.
\end{eqnarray}

\subsection{The rth estimate of $\hat{f}(x)$ and $\hat{F}(x)$}
The rth estimate of $\hat{f}(x)$ is easily obtained as follows.
\begin{eqnarray}
E(\hat{f}(x))^r&=&\int(\hat{f}(x))^r\frac{\alpha^n k^{n\alpha}}{(n-1)!}t^{-\alpha-1}[\ln(t)-n\ln(k)]^{n-1}dt\nonumber\\
&=&\int_{xk^{n-1}}^\infty\frac{(n-1)^r[\ln(t)-\ln(x)-(n-1)\ln(k)]^{r(n-2)}}{x^r[\ln(t)-n\ln(k)]^{r(n-1)}}\frac{\alpha^n k^{n\alpha}}{(n-1)!}t^{-\alpha-1}[\ln(t)-n\ln(k)]^{n-1}dt\nonumber\\
&=&\frac{(n-1)^r\alpha^n k^{\alpha n}}{x^r(n-1)!}\int_{xk^{n-1}}^\infty\frac{[\ln(t)-\ln(x)-(n-1)\ln(k)]^{r(n-2)}}{[\ln(t)-n\ln(k)]^{(r-1)(n-1)}}t^{-\alpha-1}dt.\nonumber
\end{eqnarray}
Let $z=\ln(t)-n\ln(k)\Rightarrow dz=\frac{1}{t}dt$, then
\begin{eqnarray}
E(\hat{f}(x))^r&=&\frac{(n-1)^r\alpha^n k^{\alpha n}}{x^r(n-1)!}\int_{\ln\left(\frac{x}{k}\right)}^\infty \frac{[z-\ln(x)+\ln(k)]^{r(n-2)}k^{-\alpha n}e^{-\alpha z}}{z^{(r-1)(n-1)}}dz\nonumber\\
&=&\frac{(n-1)^r\alpha^n}{x^r(n-1)!}\int_{\ln\left(\frac{x}{k}\right)}^\infty \frac{\sum_{j=0}^{r(n-2)}C(r(n-2),j)\left(-\ln\left(\frac{x}k\right)\right)^jz^{r(n-2)-j}e^{-\alpha z}}{z^{(r-1)(n-1)}}dz\nonumber\\
&=&\frac{(n-1)^r\alpha^n}{x^r(n-1)!}\sum_{j=0}^{r(n-2)}C(r(n-2),j)\left(-\ln\left(\frac{x}k\right)\right)^j
\int_{\ln\left(\frac{x}{k}\right)}^\infty z^{n-r-j-1}e^{-\alpha z}dz\nonumber\\
&=&\frac{(n-1)^r\alpha^n}{x^r(n-1)!}\sum_{j=0}^{r(n-2)}C(r(n-2),j)\left(-\ln\left(\frac{x}k\right)\right)^j
\nonumber\\&\times&
\frac{\Gamma(n-r-j)}{\alpha^{n-r-j}}\sum_{i=0}^{n-r-j-1}\frac{\exp\left(-\alpha\ln\left(\frac{x}k\right)\right)
\left(\alpha\ln\left(\frac{x}k\right)\right)^i}{i!}\nonumber
\end{eqnarray}
Therefore
\begin{eqnarray}
E(\hat{f}(x))^r=\frac{(n-1)^r\alpha^r k^\alpha}{x^{\alpha+r}(n-1)!} \sum_{j=0}^{n-r-1}C(r(n-2),j)\alpha^j\left[-\ln\left(\frac{x}{k}\right)\right]^j\Gamma(n-r-j)
\sum_{i=0}^{n-r-j-1}\frac{\alpha^i\left(\ln\left(\frac{x}k\right)\right)^i}{i!}.
\end{eqnarray}
Also, the rth estimate of $\hat{F}(x)$ is similarly obtained as follows.
\begin{eqnarray}
E(\hat{F}(x))^r&=&\int(\hat{F}(x))^rh^*(t)dt\nonumber\\
&=&\int_{xk^{n-1}}^\infty\left[1-\frac{[\ln(t)-\ln(x)-(n-1)\ln(k)]^{n-1}}{[\ln(t)-n\ln(k)]^{n-1}}\right]^rh^*(t)dt
+\int_{k^n}^{xk^{n-1}}1^r\times h^*(t)dt\nonumber\\
&=&\int_{xk^{n-1}}^\infty\sum_{j=0}^{r}C(r,j)(-1)^j\frac{[\ln(t)-\ln(x)-(n-1)\ln(k)]^{j(n-1)}}{[\ln(t)-n\ln(k)]^{j(n-1)}}
\nonumber\\&\times&\frac{\alpha^n k^{n\alpha}}{(n-1)!}t^{-\alpha-1}[\ln(t)-n\ln(k)]^{n-1}dt
+\int_{k^n}^{xk^{n-1}}\frac{\alpha^n k^{n\alpha}}{(n-1)!}t^{-\alpha-1}[\ln(t)-n\ln(k)]^{n-1}dt\nonumber\\
&=&\frac{\alpha^n k^{n\alpha}}{(n-1)!}\sum_{j=0}^{r}C(r,j)(-1)^j\int_{xk^{n-1}}^\infty
\frac{[\ln(t)-\ln(x)-(n-1)\ln(k)]^{j(n-1)}}{[\ln(t)-n\ln(k)]^{(j-1)(n-1)}}t^{-\alpha-1}dt\nonumber\\
&+&\frac{\alpha^n k^{n\alpha}}{(n-1)!}\int_{k^n}^{xk^{n-1}}t^{-\alpha-1}[\ln(t)-n\ln(k)]^{n-1}dt.\nonumber
\end{eqnarray}
Let $z=\ln(t)-n\ln(k)\Rightarrow dz=\frac{1}{t}dt$, then
\begin{eqnarray}
E(\hat{F}(x))^r&=&\frac{\alpha^n k^{n\alpha}}{(n-1)!}\sum_{j=0}^{r}C(r,j)(-1)^j\int_{\ln(\frac{x}k)}^\infty
\frac{[z-\ln(x)+\ln(k)]^{j(n-1)}}{z^{(j-1)(n-1)}}k^{-\alpha n}e^{-\alpha z}dz\nonumber\\
&+&\frac{\alpha^n k^{n\alpha}}{(n-1)!}\int_{0}^{\ln(\frac{x}k)}z^{n-1}k^{-\alpha n}e^{-\alpha z}dz\nonumber\\
&=&\frac{\alpha^n}{(n-1)!}\sum_{j=0}^{r}C(r,j)(-1)^j\int_{\ln(\frac{x}k)}^\infty
\frac{\sum_{i=0}^{j(n-1)}C(j(n-1),i)(-\ln(\frac{x}k))^i z^{j(n-1)-i}e^{-\alpha z}}{z^{(j-1)(n-1)}}dz\nonumber\\
&+&\frac{\alpha^n}{(n-1)!}\frac{\Gamma(n)}{\alpha^n}\left[1-\left(\frac{x}k\right)^{-\alpha}
\sum_{i=0}^{n-1}\frac{(\alpha\ln(\frac{x}k))^i}{i!}\right]\nonumber\\
&=&\frac{\alpha^n}{(n-1)!}\sum_{j=0}^{r}C(r,j)(-1)^j\sum_{i=0}^{j(n-1)}C(j(n-1),i)\left(-\ln\left(\frac{x}k\right)\right)^i
\int_{\ln(\frac{x}k)}^\infty
z^{n-i-1}e^{-\alpha z}dz\nonumber\\
&+&\left[1-\left(\frac{x}k\right)^{-\alpha}\sum_{i=0}^{n-1}\frac{(\alpha\ln(\frac{x}k))^i}{i!}\right]\nonumber\\
&=&\frac{\alpha^n}{(n-1)!}\sum_{j=0}^{r}C(r,j)(-1)^j\sum_{i=0}^{j(n-1)}C(j(n-1),i)\left(-\ln\left(\frac{x}k\right)\right)^i
\frac{\Gamma(n-i)}{\alpha^{n-i}}\nonumber\\
&\times&\sum_{l=0}^{n-i-1}\frac{\exp(-\alpha\ln(\frac{x}k))(\alpha\ln(\frac{x}k))^l}{l!}
+\left(\frac{x}k\right)^{-\alpha}
\left[\left(\frac{x}k\right)^{\alpha}-\sum_{i=0}^{n-1}\frac{(\alpha\ln(\frac{x}k))^i}{i!}\right].\nonumber
\end{eqnarray}
Then
\begin{eqnarray}
E(\hat{F}(x))^r&=&\frac{k^\alpha}{(n-1)!x^\alpha}\sum_{j=0}^{r}C(r,j)(-1)^j\sum_{i=0}^{j(n-1)}C(j(n-1),i)
\left(-\ln\left(\frac{x}k\right)\right)^i\alpha^i\Gamma(n-i)\nonumber\\
&\times&\sum_{l=0}^{n-i-1}\frac{(\alpha\ln(\frac{x}k))^l}{l!}+\left(\frac{k}x\right)^{\alpha}
\left[\left(\frac{x}k\right)^{\alpha}-\sum_{i=0}^{n-1}\frac{(\alpha\ln(\frac{x}k))^i}{i!}\right].
\end{eqnarray}

\section{R code}
The R code to compare the bias and MSE of the estimators is as follows.
\begin{verbatim}
sim=function(t,n,k,alpha,r)
{
sfh<-0
sFh<-0
sft<-0
sFt<-0
for(l in 1:t){
x<-array(, c(1,n))
for (i in 1:n)          {
u<-runif(1,0,1)
x[i]<-k*(1-u)^(-1/alpha)}
alphah<-n/sum(log(x)-log(k))
fx<-alpha*k^alpha/x[1]^(alpha+1)
intB<-function(z){ z^(-n)*exp(-alpha*n/z)*(k/x[1])^z}
B<-integrate(intB,lower=0,upper=Inf)$value
Eftildx<-(alpha*n)^n/factorial(n-1)/x[1]*B
   intB1<-function(z){ z^(-n+1)*exp(-alpha*n/z)*(k/x[1])^(2*z)}
   B1<-integrate(intB1,lower=0,upper=Inf)$value
   Eftildxs2<-(alpha*n)^n/factorial(n-1)/x[1]^2*B1
MSEftildx<-Eftildxs2-2*fx*Eftildx+fx^2
intA<-function(z){ z^(2*n-4)*exp(-alpha*z)/(z+log(x[1])-log(k))^(n-1)}
A<-integrate(intA,lower=0,upper=Inf)$value
MSEfhx<-(n-1)*alpha^n*k^alpha/x[1]^(alpha+2)/factorial(n-2)*A-alpha^2*
k^(2*alpha)/x[1]^(2*alpha+2)

Fx<-1-(k/x[1])^alpha
intB2<-function(w){ w^(-n-1)*exp(w*log(k/x[1]))*exp(-alpha*n/w)}
B2<-integrate(intB2,lower=0,upper=Inf)$value
EFtildx<-1-(alpha*n)^n/factorial(n-1)*B2
      intB3<-function(w){ w^(-n-1)*exp(2*w*log(k/x[1]))*exp(-alpha*n/w)}
      B3<-integrate(intB3,lower=0,upper=Inf)$value
      EFtildxs2<-1-2*(alpha*n)^n/factorial(n-1)*B2+(alpha*n)^n/factorial(n-1)*B3
MSEFtildx<-EFtildxs2-2*Fx*EFtildx+Fx^2
intA1<-function(w){ (w-log(x[1]/k))^(2*n-2)*exp(-alpha*w)*w^(-n+1)}
A1<-integrate(intA1,lower=log(x[1]/k),upper=Inf)$value
MSEFhx<-1-2*(k/x[1])^alpha+alpha^n/factorial(n-1)*A1-(1-(k/x[1])^alpha)^2

sfh<-sfh+MSEfhx
sft<-sft+MSEftildx
sFh<-sFh+MSEFhx
sFt<-sFt+MSEFtildx
}
mMSEfhx<-sfh/t
mMSEftildx<-sft/t
mMSEFhx<-sFh/t
mMSEFtildx<-sFt/t
return(c(mMSEfhx,mMSEftildx,mMSEFhx,mMSEFtildx))
}
sim(10,5,1,5,1)

sim1=function(t,k,alpha,r){
i<-seq(3,35,1)
for (j in i){
sim(t,j,k,alpha,r)}}
sim1(10,1,5,1)
\end{verbatim}

\section{Tables}
~~~~In order to get the idea of efficiency between the two type of estimation i.e MLE and UMVUE. We have generated a sample of size 4(1)15(5)100 from the Pareto distribution with $\alpha$=0.5(0.5)2 and $k$=0.5(0.5)2. We have given Tables based on one thousand independent replication of each experiments.\\
Table 1. shows the bias and MSE of the estimators of the pdf and bias and MSE of the estimators of cdf are shown in Tables 2. The value in the bracket is for the MSE in each tables.
From the Tables, it has been seen that MLE of pdf and cdf are more efficient than UMVUEs.\\
One should note that UMVUE of $\alpha$ is better than MLE of $\alpha$.

\begin{center}
\textbf{Table 1.} MSE of $\hat{f}(x)$ and $\tilde{f}(x)$ for different values of $\alpha$ and $k$ respect to $n$
\tiny
\begin{tabular}{|l|l|l|l|l|l|l|}
\hline
$n$&~~~~$\alpha=0.5$&~~~~~$\alpha=1$&~$\alpha=1.5$&~~$\alpha=2$&~~$\alpha=0.5$&~~$\alpha=2$\\
&~~~~$k=0.5$&~~~~~$k=1$&~$k=1.5$&~~$k=2$&~~$k=2$&~~$k=0.5$\\
\hline
4&.4723690000&.4755806957&.476511&.476705&.029389&7.658800\\
&(.4551780000)&(.4572692552)&(.457605)&(.457264)&(.028230)&(7.374350)\\
\hline
5&.3172310000&.3207780187&.322473&.321924&.019827&5.157680\\
&(.2725650000)&(.2757554306)&(.277487)&(.276551)&(.017036)&(4.434020)\\
\hline
6&.2383780000&.2404067792&.241698&.242471&.014881&3.885290\\
&(.1881800000)&(.1896539273)&(.190811)&(.191536)&(.011743)&(3.071030)\\
\hline
7&.1914040000&.1929655509&.193791&.194454&.011901&3.105110\\
&(.1421410000)&(.1432299232)&(.143884)&(.144461)&(.008828)&(2.305270)\\
\hline
8&.1599510000&.1602595912&.161806&.161604&.009946&2.588020\\
&(.1133700000)&(.1133783291)&(.114652)&(.114415)&(.007042)&(1.832800)\\
\hline
9&.1372420000&.1380770369&.138780&.139051&.008601&2.213210\\
&(.0937400000)&(.0942468652)&(.094764)&(.094956)&(.005877)&(1.509400)\\
\hline
10&.1195130000&.1209664617&.121456&.121541&.007514&1.944830\\
&(.0791650000)&(.0801522609)&(.080488)&(.080528)&(.004982)&(1.288590)\\
\hline
11&.1063690000&.1075933345&.107972&.108015&.006661&1.728000\\
&(.0687560000)&(.0695600949)&(.069808)&(.069818)&(.004307)&(1.116900)\\
\hline
12&.0957350000&.0966263202&.096937&.097635&.005985&1.561640\\
&(.0606160000)&(.0611705833)&(.061363)&(.061865)&(.003789)&(.989450)\\
\hline
13&.0873630000&.0876848730&.088438&.088472&.005442&1.418000\\
&(.0543730000)&(.0545292384)&(.055038)&(.055049)&(.003386)&(.882550)\\
\hline
14&.0802730000&.0808533231&.081013&.081069&.005006&1.299260\\
&(.0492200000)&(.0495643965)&(.049654)&(.049681)&(.003069)&(.796410)\\
\hline
15&.0735700000&.0746156473&.074685&.074919&.004613&1.197330\\
&(.0444890000)&(.0451395324)&(.045168)&(.045318)&(.002790)&(.724140)\\
\hline
20&.0532991000&.0538278427&.053924&.054039&.003330&.865418\\
&(.0307735000)&(.0310780487)&(.031129)&(.031196)&(.001923)&(.499644)\\
\hline
25&.0418259000&.0420976064&.042261&.042273&.002612&.676949\\
&(.0234799000)&(.0236273808)&(.023720)&(.023724)&(.001466)&(.379933)\\
\hline
30&.0341756000&.0346115111&.034740&.034775&.002141&.556918\\
&(.0188200000)&(.0190619643)&(.019133)&(.019152)&(.001179)&(.306728)\\
\hline
35&.0290451000&.0294767054&.029444&.029448&.001821&.473159\\
&(.0157795000)&(.0160172951)&(.015996)&(.015996)&(.000989)&(.257090)\\
\hline
40&.0252562000&.0254701187&.025650&.025593&.001581&.411040\\
&(.0135817000)&(.0136957411)&(.013795)&(.013761)&(.000850)&(.221060)\\
\hline
45&.0224046000&.0225300956&.022630&.022658&.001395&.362473\\
&(.0119535000)&(.0120188221)&(.012073)&(.012088)&(.000744)&(.193366)\\
\hline
50&.0198892000&.0202331230&.020334&.020295&.001253&.324570\\
&(.0105411000)&(.0107253823)&(.010780)&(.010758)&(.000664)&(.172036)\\
\hline
55&.0181373000&.0183714478&.018323&.018382&.001130&.294111\\
&(.0095638000)&(.0096881397)&(.009661)&(.009692)&(.000596)&(.155077)\\
\hline
60&.0165965000&.0167227711&.016761&.016767&.001036&.269022\\
&(.0087132000)&(.0087790315)&(.008799)&(.008802)&(.000544)&(.141232)\\
\hline
65&.0152848000&.0153716098&.015452&.015489&.000954&.247479\\
&(.0079948000)&(.0080394109)&(.008082)&(.008102)&(.000499)&(.129436)\\
\hline
70&.0141436000&.0142593414&.014292&.014337&.000880&.228735\\
&(.0073741000)&(.0074341631)&(.007451)&(.007474)&(.000459)&(.119243)\\
\hline
75&.0131255000&.0132738528&.013345&.013349&.000820&.213970\\
&(.0068239000)&(.0069011346)&(.006938)&(.006940)&(.000427)&(.111250)\\
\hline
80&.0123438000&.0124114845&.012482&.012486&.000766&.199643\\
&(.0064025000)&(.0064370354)&(.006474)&(.006476)&(.000397)&(.103542)\\
\hline
85&.0115424000&.0117098555&.011725&.011738&.000722&.188222\\
&(.0059735000)&(.0060606327)&(.006068)&(.006075)&(.000374)&(.097417)\\
\hline
90&.0109372000&.0110042351&.011088&.011080&.000680&.177858\\
&(.0056499000)&(.0056841841)&(.005728)&(.005724)&(.000351)&(.091881)\\
\hline
95&.0102967000&.0104005835&.010482&.010493&.000644&.167683\\
&(.0053096000)&(.0053631148)&(.005405)&(.005411)&(.000332)&(.086471)\\
\hline
100&.0098069000&.0098648715&.009939&.009971&.000615&.159213\\
&(.0050495000)&(.0050790332)&(.005117)&(.005134)&(.000316)&(.081978)\\
\hline
\end{tabular}

The figures in the bracket refers to the MSE of MLE of $f(x)$ ($\tilde{f}(x)$) and without bracket refers to the MSE of UMVUE of $f(x)$ ($\hat{f}(x)$)
\end{center}
\newpage
\begin{center}
\textbf{Table 2.} MSE of $\hat{F}(x)$ and $\tilde{F}(x)$ for different values of $\alpha$ and $k$ respect to $n$
\tiny
\begin{tabular}{|l|l|l|l|l|l|l|}
\hline
$n$&~~~~$\alpha=0.5$&~~~~~$\alpha=1$&~$\alpha=1.5$&~~$\alpha=2$&~~$\alpha=0.5$&~~$\alpha=2$\\
&~~~~$k=0.5$&~~~~~$k=1$&~$k=1.5$&~~$k=2$&~~$k=2$&~~$k=0.5$\\
\hline
4&.1333908446&.1368410441&.038435&.115906&.121205&.157459\\
&(.0014763889)&(.0012562228)&(.001100)&(.000794)&(.001880)&(.002199)\\
\hline
5&.1118961062&.1799791911&.172372&.198062&.137497&.124812\\
&(.0072949912)&(.0356532489)&(.029683)&(.055201)&(.013072)&(.009786)\\
\hline
6&.1408169615&.1737093067&.101804&.188658&.143735&.134947\\
&(.0037596610)&(.0032023158)&(.001684)&(.005213)&(.003836)&(.005169)\\
\hline
7&.2191326212&.2065645792&.195758&.250334&.126943&.160016\\
&(.0229514291)&(.0165453192)&(.013057)&(.067911)&(.004378)&(.007159)\\
\hline
8&.2007151499&.2012409768&.195431&.219443&.197351&.203045\\
&(.0034259269)&(.0032392638)&(.004168)&(.000621)&(.002333)&(.007496)\\
\hline
9&.2576817407&.2512987261&.147364&.135785&.124131&.269592\\
&(.0275929634)&(.0205299145)&(.004077)&(.003471)&(.002910)&(.060899)\\
\hline
10&.1492978531&.1837933180&.150179&.220440&.134771&.213361\\
&(.0023002929)&(.0030831665)&(.003979)&(.005504)&(.001508)&(.007332)\\
\hline
11&.2623413820&.2259891179&.259195&.274008&.163600&.129209\\
&(.0139901431)&(.0076116749)&(.012447)&(.046732)&(.003893)&(.002419)\\
\hline
12&.1819439058&.1679376346&.197488&.176699&.147436&.205010\\
&(.0062366813)&(.0049937937)&(.004988)&(.004595)&(.004263)&(.005586)\\
\hline
13&.2177045578&.1674005776&.235847&.254445&.141234&.141272\\
&(.0058142184)&(.0033358805)&(.006945)&(.008379)&(.002349)&(.002350)\\
\hline
14&.2356903001&.0548103487&.227515&.193420&.188948&.136300\\
&(.0060898969)&(.0003127263)&(.006204)&(.004149)&(.006286)&(.001994)\\
\hline
15&.1076037454&.1555011068&.122969&.253257&.199792&.191199\\
&(.0011259936)&(.0024132218)&(.001482)&(.007139)&(.004118)&(.003743)\\
\hline
20&.2627256864&.2622509943&.248371&.232384&.164861&.168388\\
&(.0065292099)&(.0062066109)&(.005120)&(.004292)&(.002324)&(.002052)\\
\hline
25&.2488455377&.2598193829&.218770&.259375&.151849&.177311\\
&(.0041453670)&(.0052022850)&(.004877)&(.004883)&(.001275)&(.001795)\\
\hline
30&.2251732413&.2172670151&.257171&.211684&.113000&.118642\\
&(.0042208658)&(.0023800197)&(.004369)&(.003983)&(.000552)&(.000612)\\
\hline
35&.2458125666&.2072001232&.198559&.204366&.178507&.156018\\
&(.0028694984)&(.0018011252)&(.003254)&(.003352)&(.001266)&(.000935)\\
\hline
40&.1737078950&.2554091293&.251794&.152611&.116046&.156706\\
&(.0010316986)&(.0032353869)&(.002759)&(.002133)&(.001487)&(.000818)\\
\hline
45&.1940551896&.2429150460&.173837&.222621&.157868&.108663\\
&(.0011808928)&(.0021535385)&(.002250)&(.002867)&(.000733)&(.000325)\\
\hline
50&.1064113870&.2267825379&.179884&.228218&.125222&.151083\\
&(.0011050152)&(.0015767910)&(.000884)&(.002635)&(.000394)&(.000594)\\
\hline
55&.2299332896&.1855659085&.252088&.218652&.163238&.181364\\
&(.0024125373)&(.0008607060)&(.002056)&(.002338)&(.000639)&(.001967)\\
\hline
60&.1392457773&.1773309892&.127104&.227032&.170991&.108474\\
&(.0013392528)&(.0017724707)&(.001195)&(.002202)&(.000649)&(.000238)\\
\hline
65&.1986728058&.2322009300&.253793&.115857&.107674&.118045\\
&(.0018421105)&(.0012944937)&(.001897)&(.000989)&(.000215)&(.000262)\\
\hline
70&.2288138955&.1496013496&.251293&.108626&.134125&.123655\\
&(.0019021751)&(.0012705637)&(.001864)&(.000203)&(.000320)&(.000268)\\
\hline
75&.2524003249&.1656385542&.252058&.226321&.113991&.124311\\
&(.0017059572)&(.0013387601)&(.001537)&(.001769)&(.000209)&(.000252)\\
\hline
80&.2232959249&.1154621504&.238145&.123070&.159104&.104494\\
&(.0016486241)&(.0008181671)&(.001686)&(.000887)&(.000408)&(.000719)\\
\hline
85&.1569167081&.1116721695&.197303&.038435&.100752&.160423\\
&(.0003712149)&(.0007421171)&(.000637)&(.000164)&(.000141)&(.000959)\\
\hline
90&.2498099645&.1551047337&.154748&.084016&.080884&.071267\\
&(.0014551335)&(.0003408055)&(.000339)&(.000483)&(.000459)&(.000064)\\
\hline
95&.2525238629&.1423488455&.237976&.038663&.021378&.055160\\
&(.0013016097)&(.0002657662)&(.000956)&(.000151)&(.000005)&(.000255)\\
\hline
100&.1352354632&.0786840925&.105801&.064531&.055326&.026498\\
&(.0008108742)&(.0004021173)&(.000132)&(.000305)&(.000245)&(.000081)\\
\hline
\end{tabular}

The figures in the bracket refers to the MSE of MLE of $F(x)$ ($\tilde{F}(x)$) and without bracket refers to the MSE of UMVUE of $F(x)$ ($\hat{F}(x)$)
\end{center}

\newpage
\input{epsf}
\begin{figure}
\begin{tabular}{rr}
\epsfxsize=6.5in\epsfysize=10in \epsffile{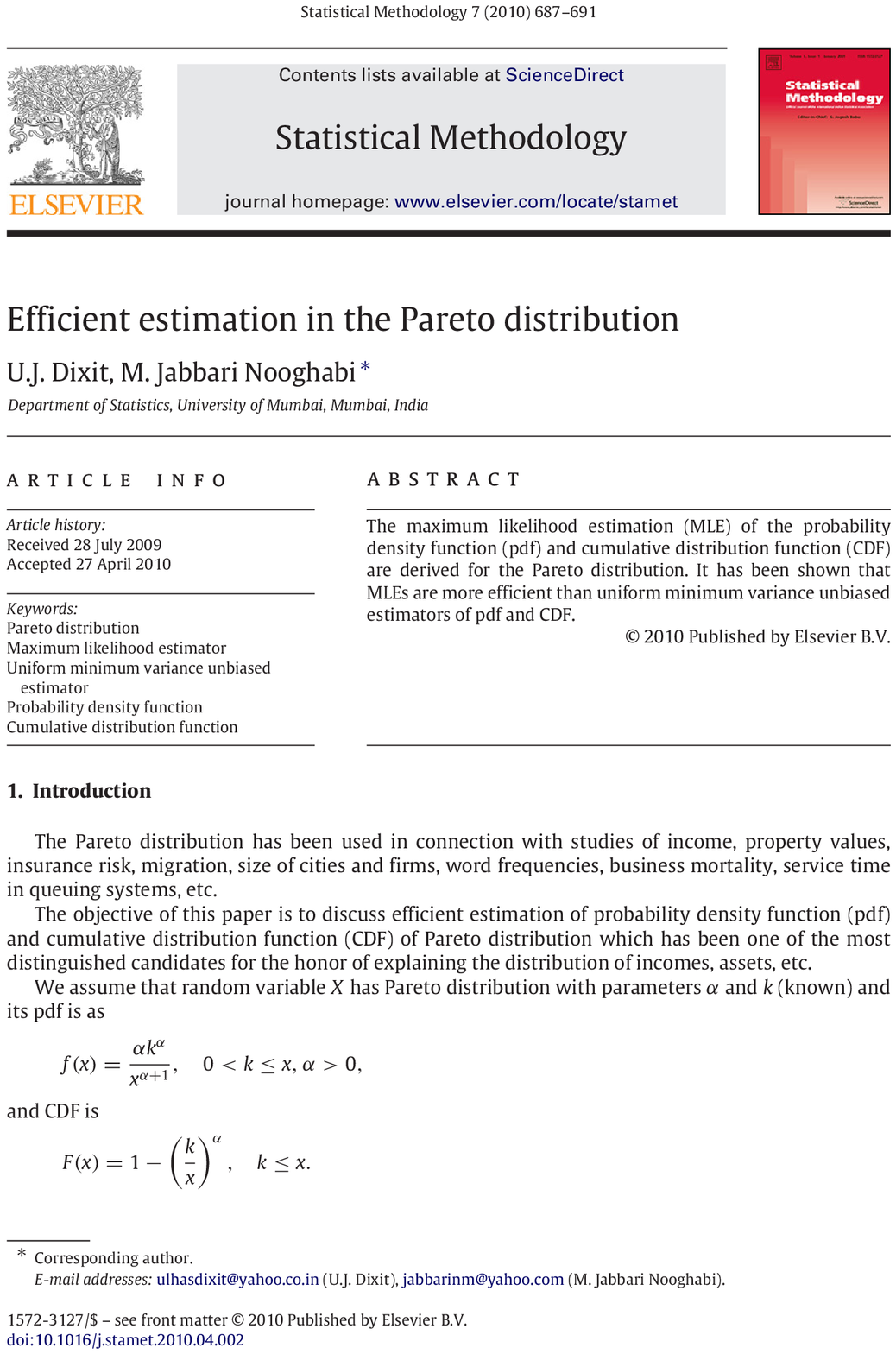}
\end{tabular}
\end{figure}
\input{epsf}
\begin{figure}
\begin{tabular}{rr}
\epsfxsize=6.5in\epsfysize=10in \epsffile{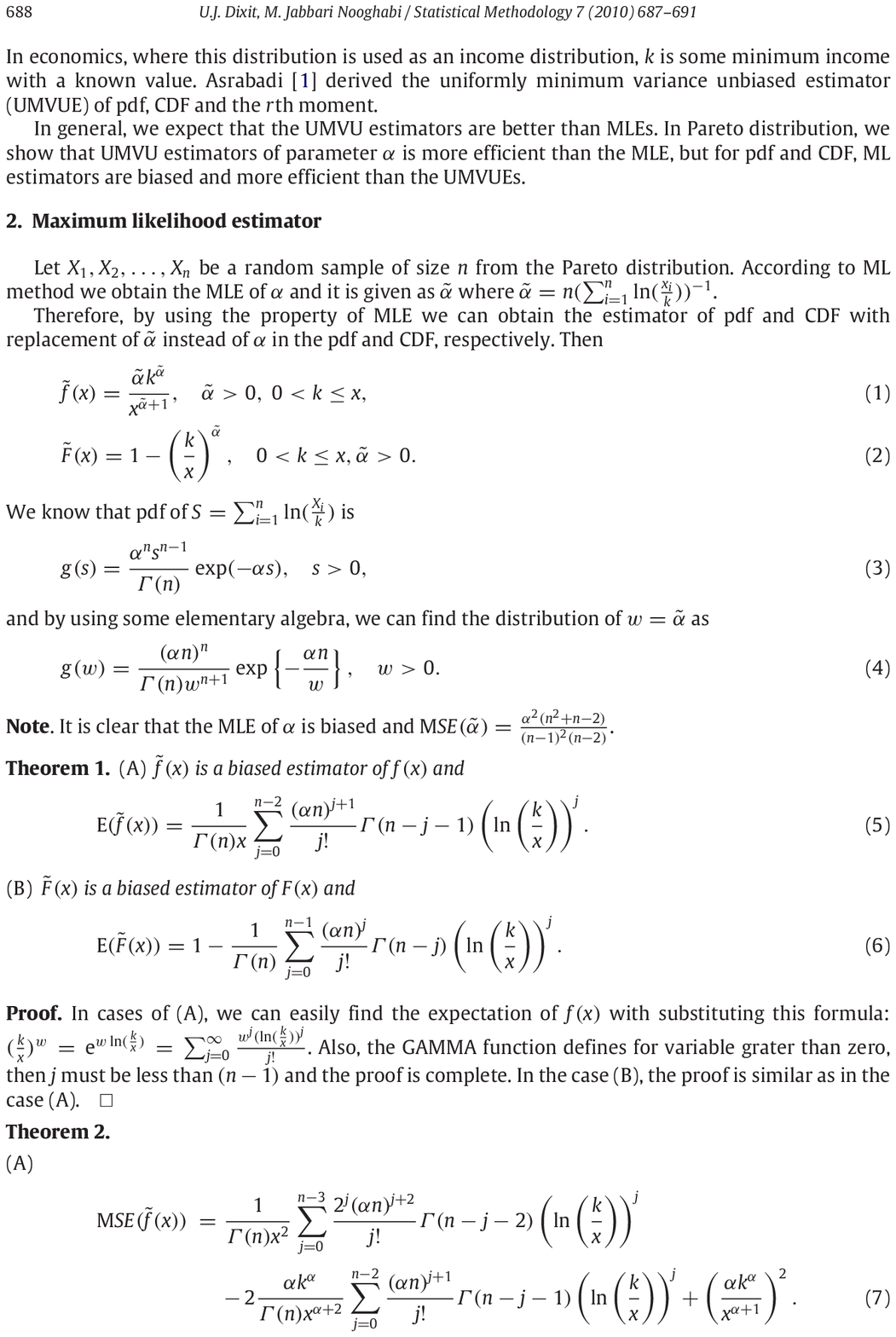}
\end{tabular}
\end{figure}
\input{epsf}
\begin{figure}
\begin{tabular}{rr}
\epsfxsize=6.5in\epsfysize=10in \epsffile{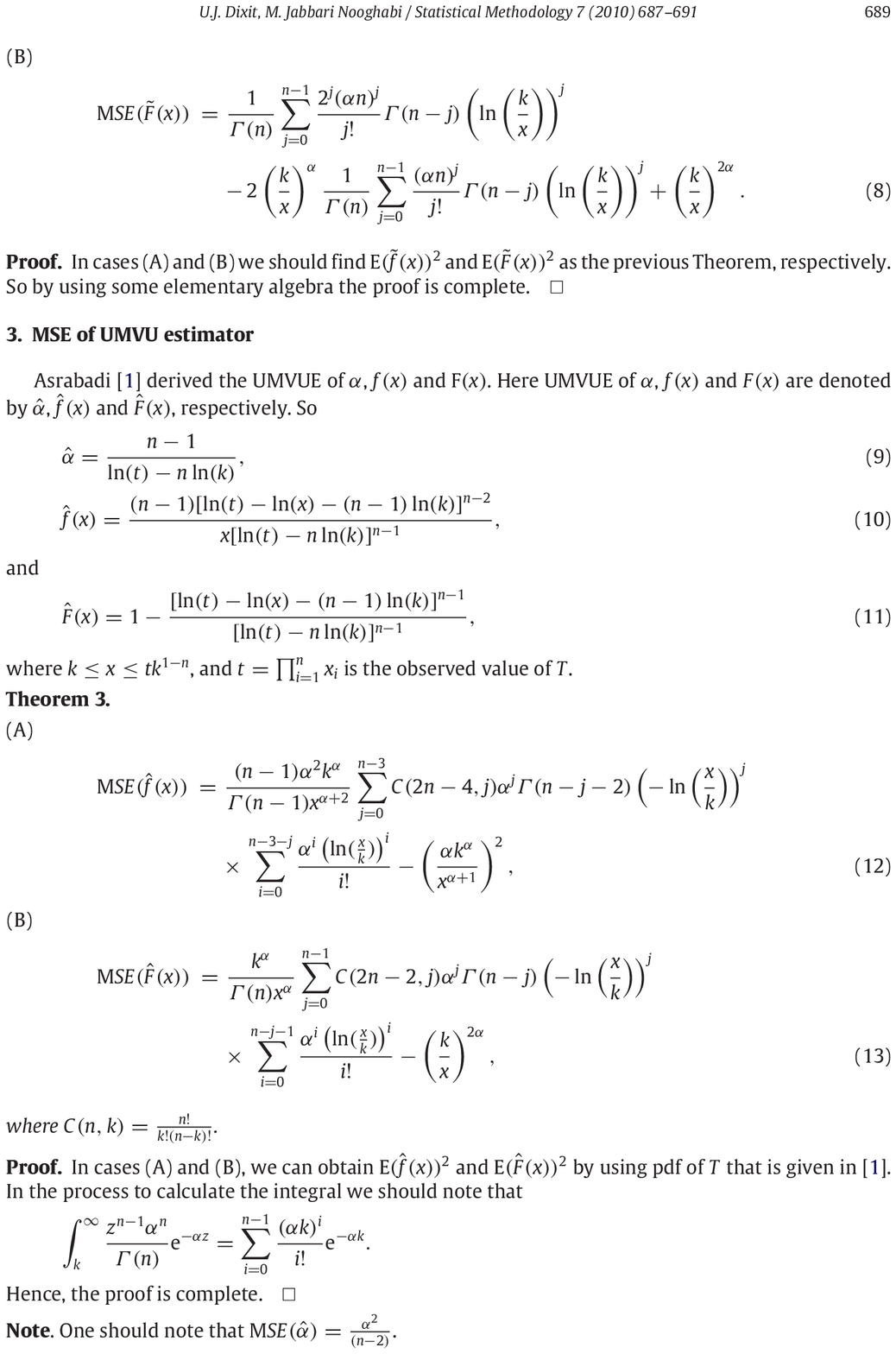}
\end{tabular}
\end{figure}
\input{epsf}
\begin{figure}
\begin{tabular}{rr}
\epsfxsize=6.5in\epsfysize=10in \epsffile{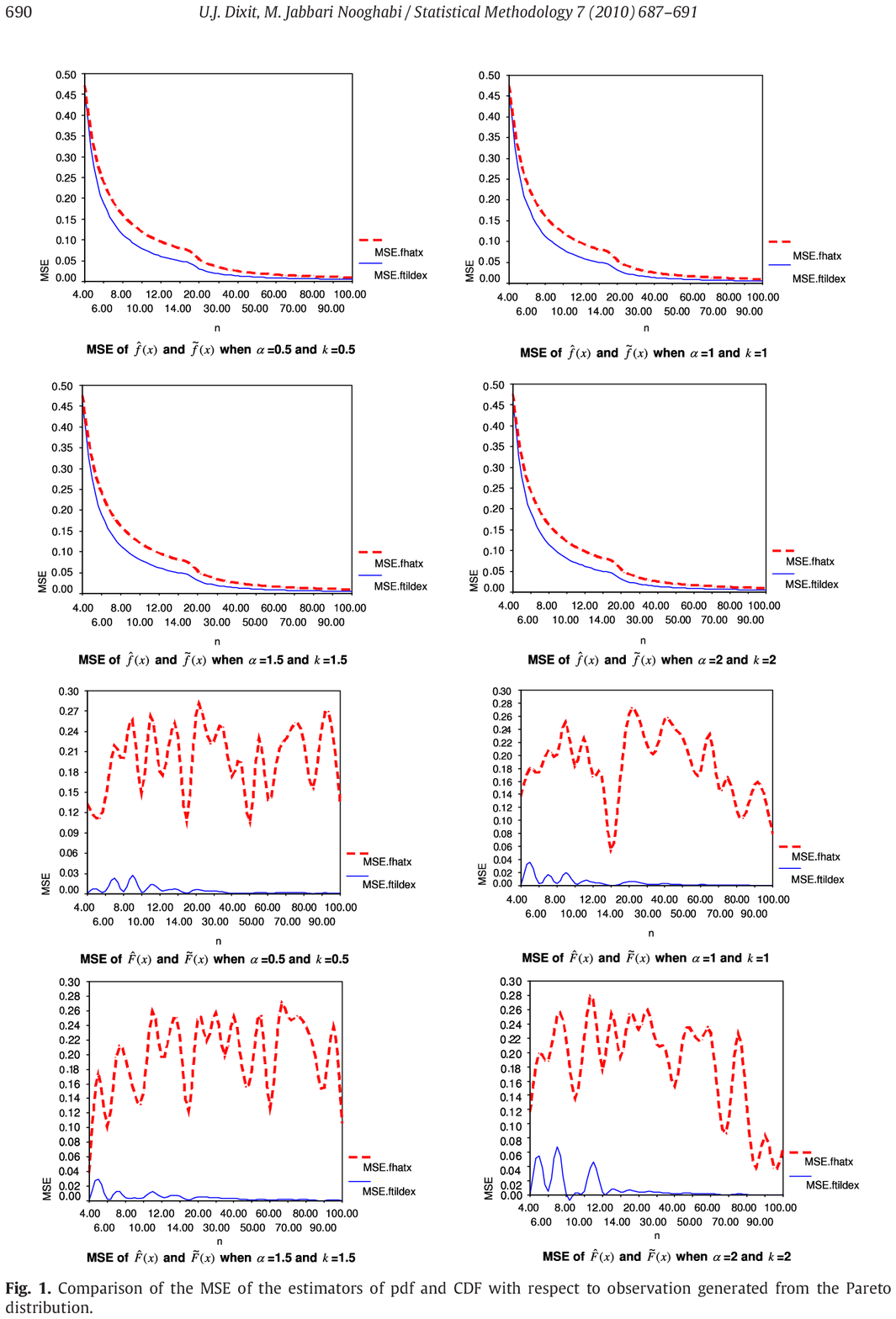}
\end{tabular}
\end{figure}
\input{epsf}
\begin{figure}
\begin{tabular}{rr}
\epsfxsize=6.5in\epsfysize=10in \epsffile{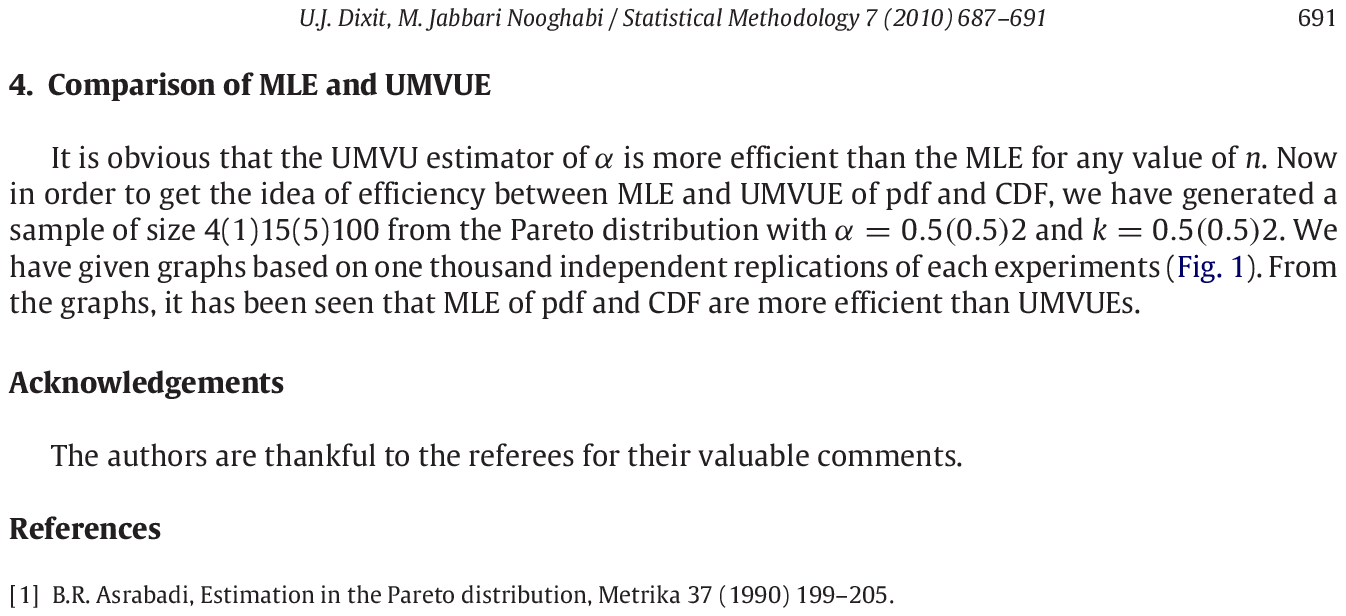}
\end{tabular}
\end{figure}
\input{epsf}
\begin{figure}
\begin{tabular}{rr}
\epsfxsize=6.5in\epsfysize=10in \epsffile{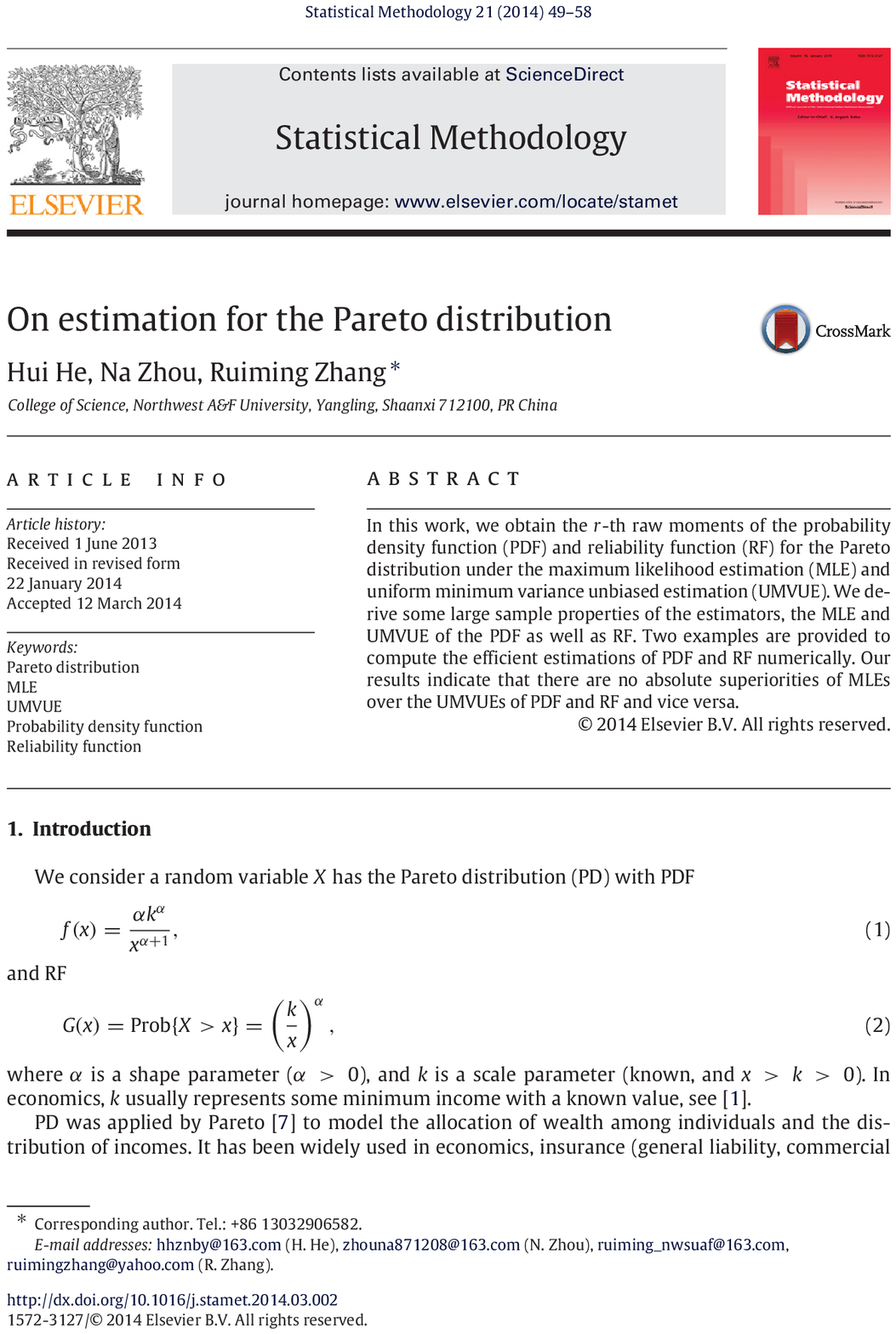}
\end{tabular}
\end{figure}
\input{epsf}
\begin{figure}
\begin{tabular}{rr}
\epsfxsize=6.5in\epsfysize=10in \epsffile{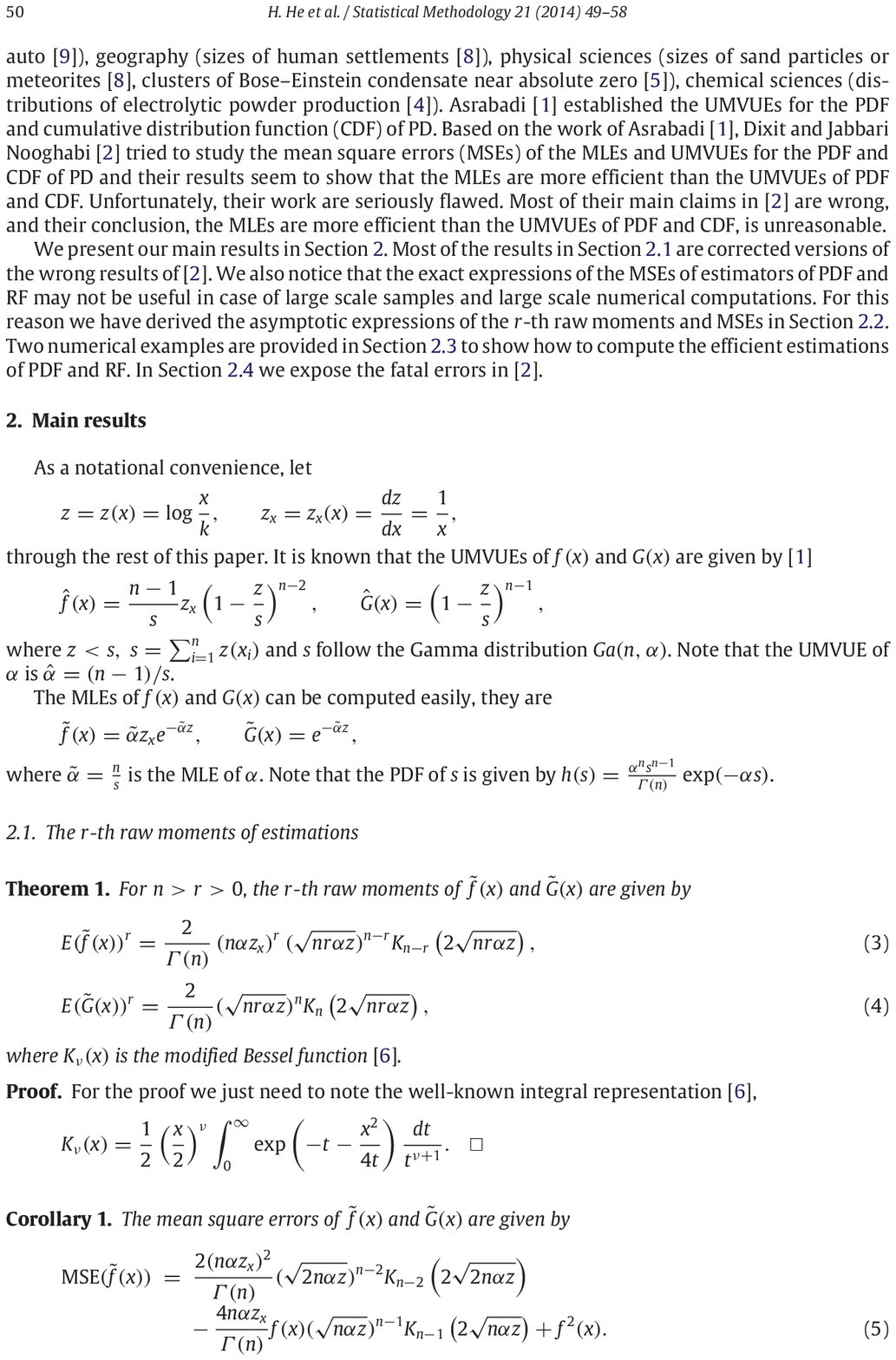}
\end{tabular}
\end{figure}
\input{epsf}
\begin{figure}
\begin{tabular}{rr}
\epsfxsize=6.5in\epsfysize=10in \epsffile{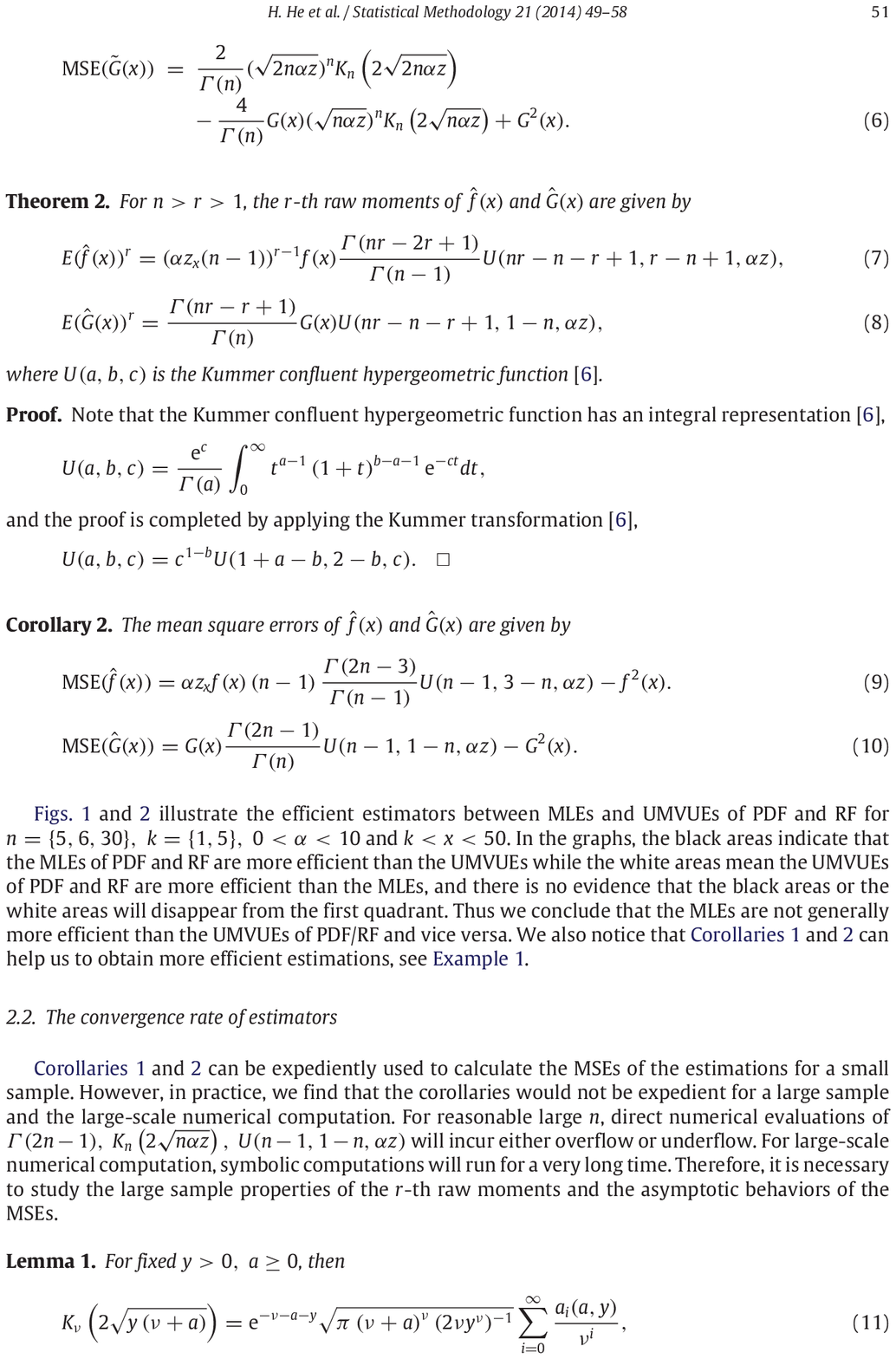}
\end{tabular}
\end{figure}
\input{epsf}
\begin{figure}
\begin{tabular}{rr}
\epsfxsize=6.5in\epsfysize=10in \epsffile{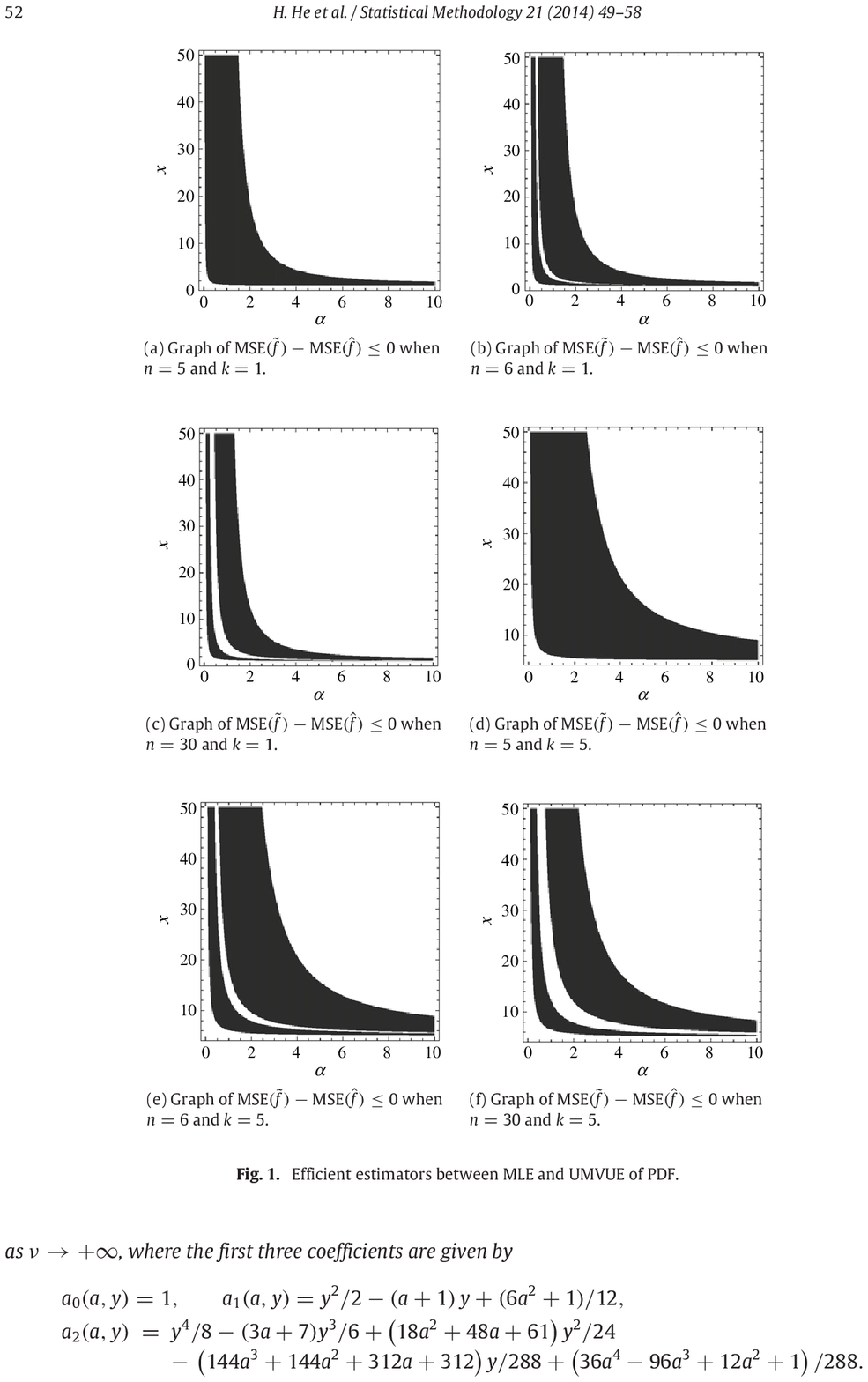}
\end{tabular}
\end{figure}
\input{epsf}
\begin{figure}
\begin{tabular}{rr}
\epsfxsize=6.5in\epsfysize=10in \epsffile{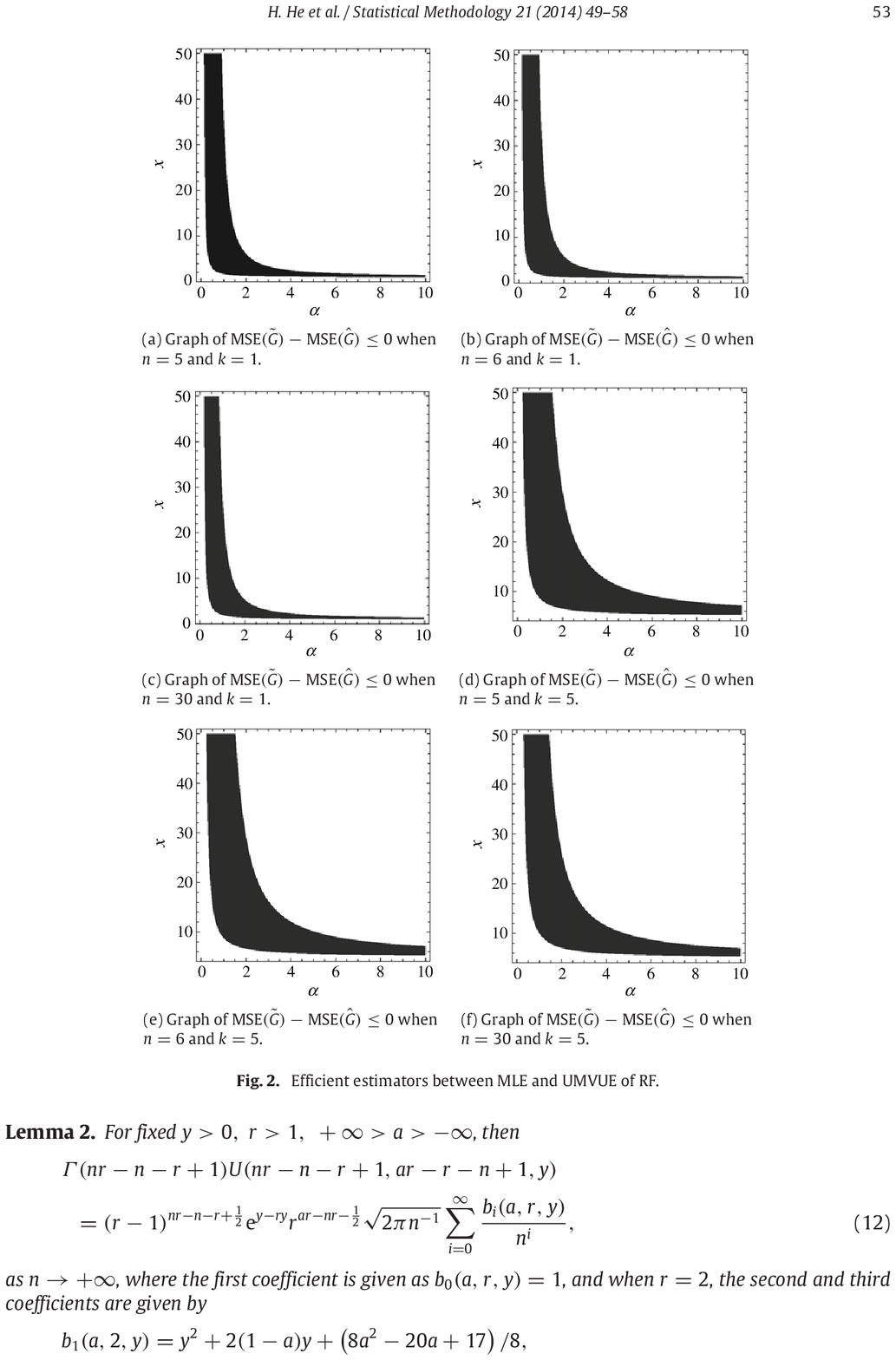}
\end{tabular}
\end{figure}
\input{epsf}
\begin{figure}
\begin{tabular}{rr}
\epsfxsize=6.5in\epsfysize=10in \epsffile{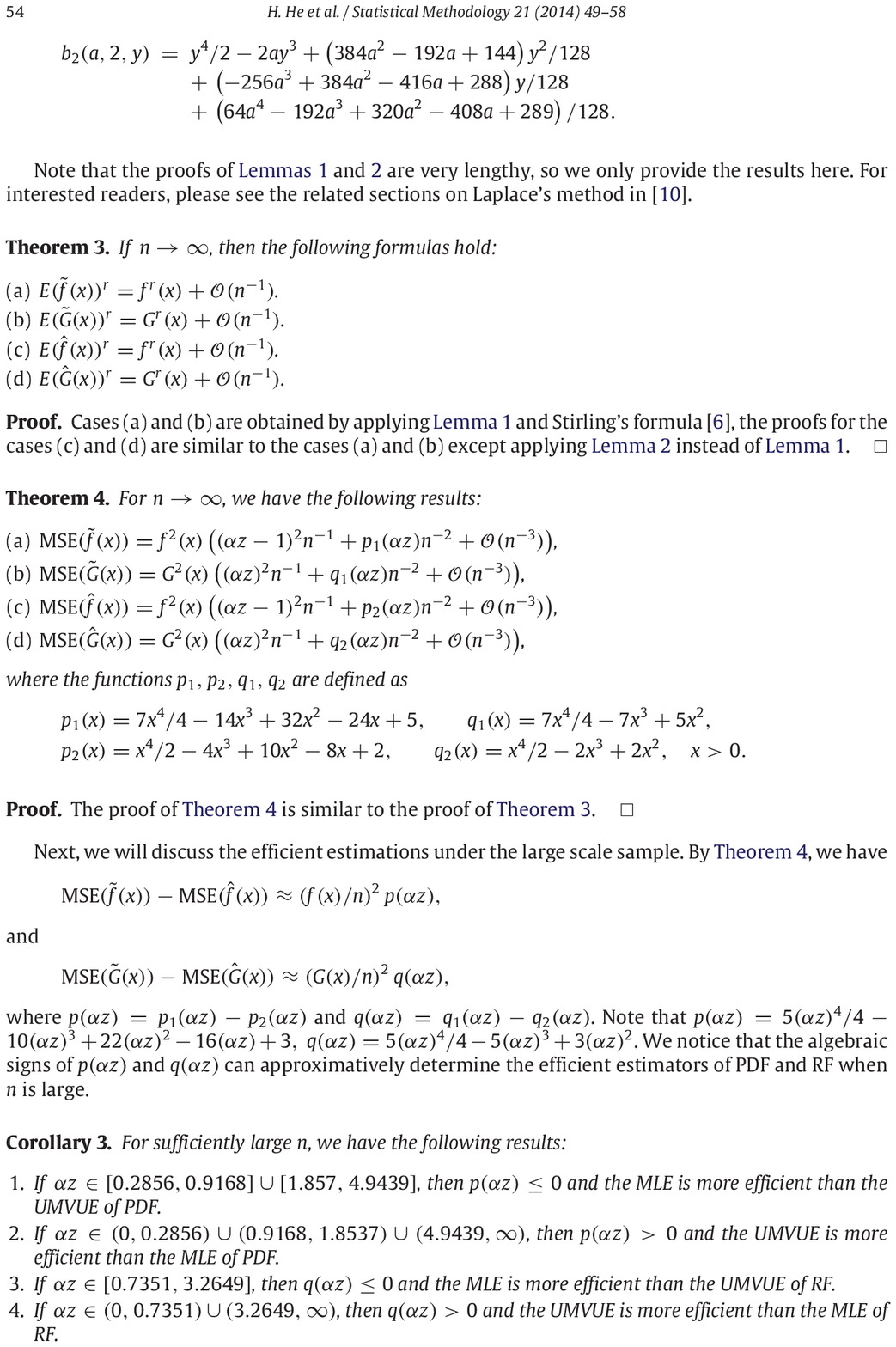}
\end{tabular}
\end{figure}
\input{epsf}
\begin{figure}
\begin{tabular}{rr}
\epsfxsize=6.5in\epsfysize=10in \epsffile{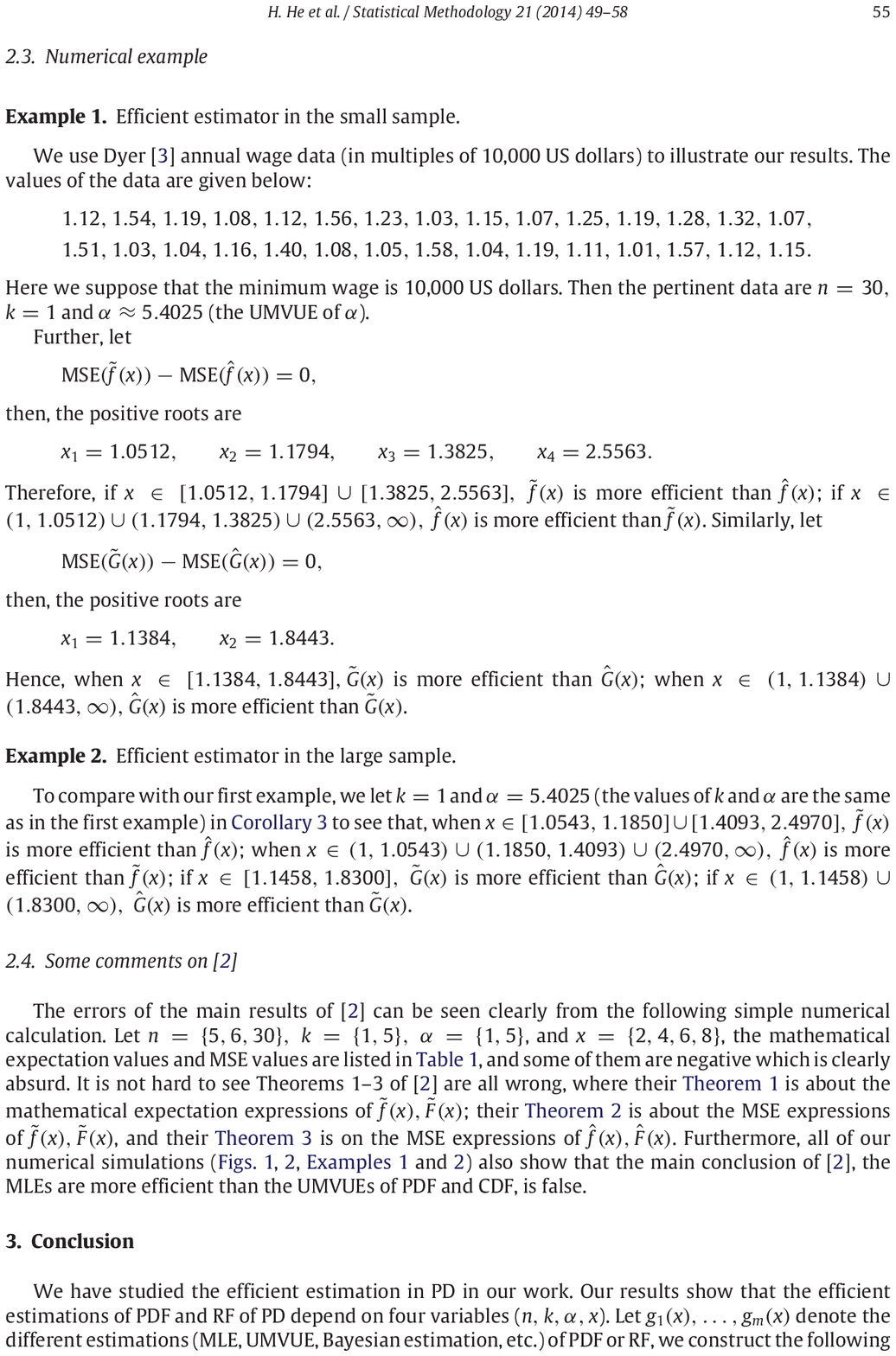}
\end{tabular}
\end{figure}
\input{epsf}
\begin{figure}
\begin{tabular}{rr}
\epsfxsize=6.5in\epsfysize=10in \epsffile{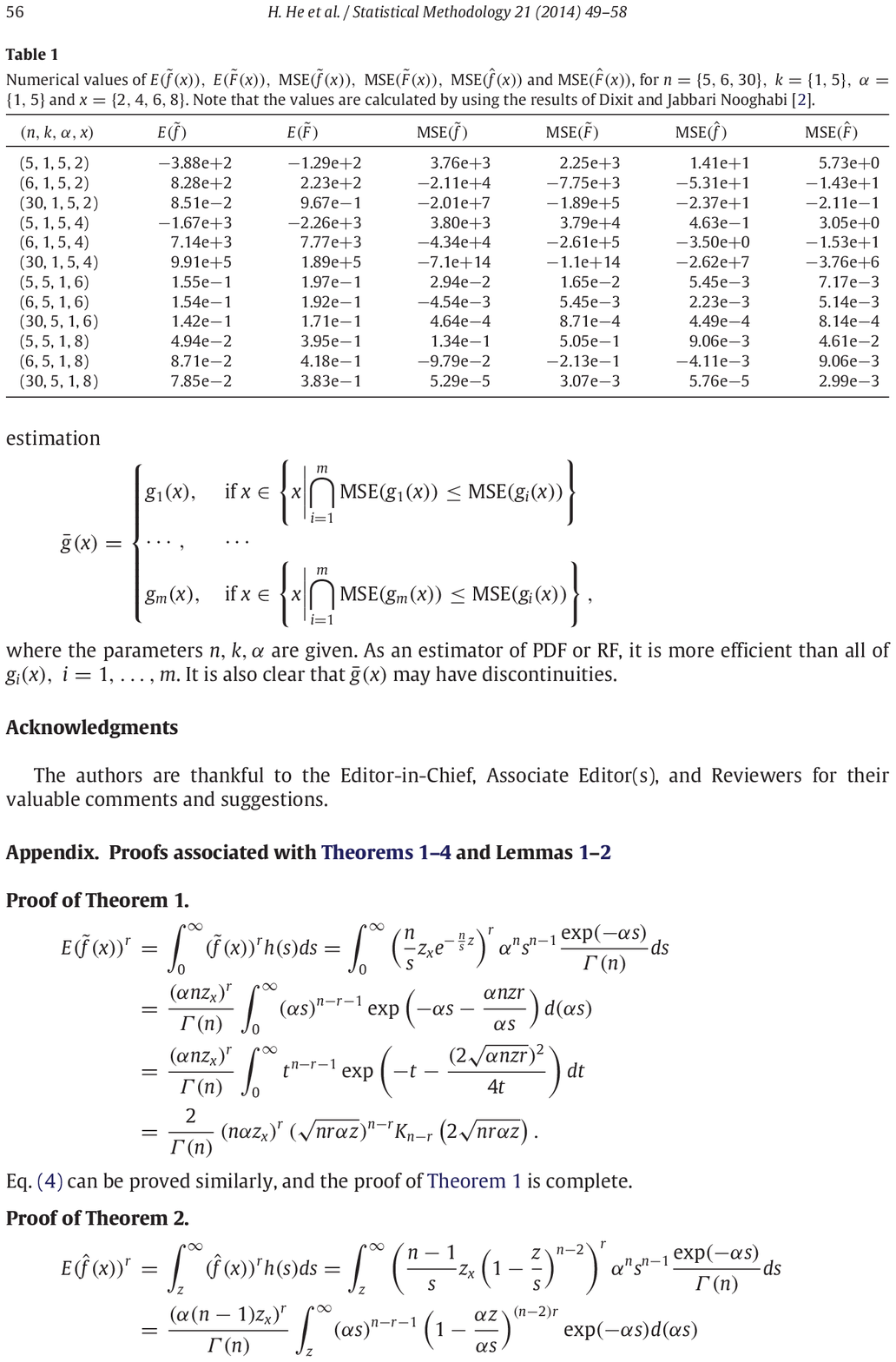}
\end{tabular}
\end{figure}
\input{epsf}
\begin{figure}
\begin{tabular}{rr}
\epsfxsize=6.5in\epsfysize=10in \epsffile{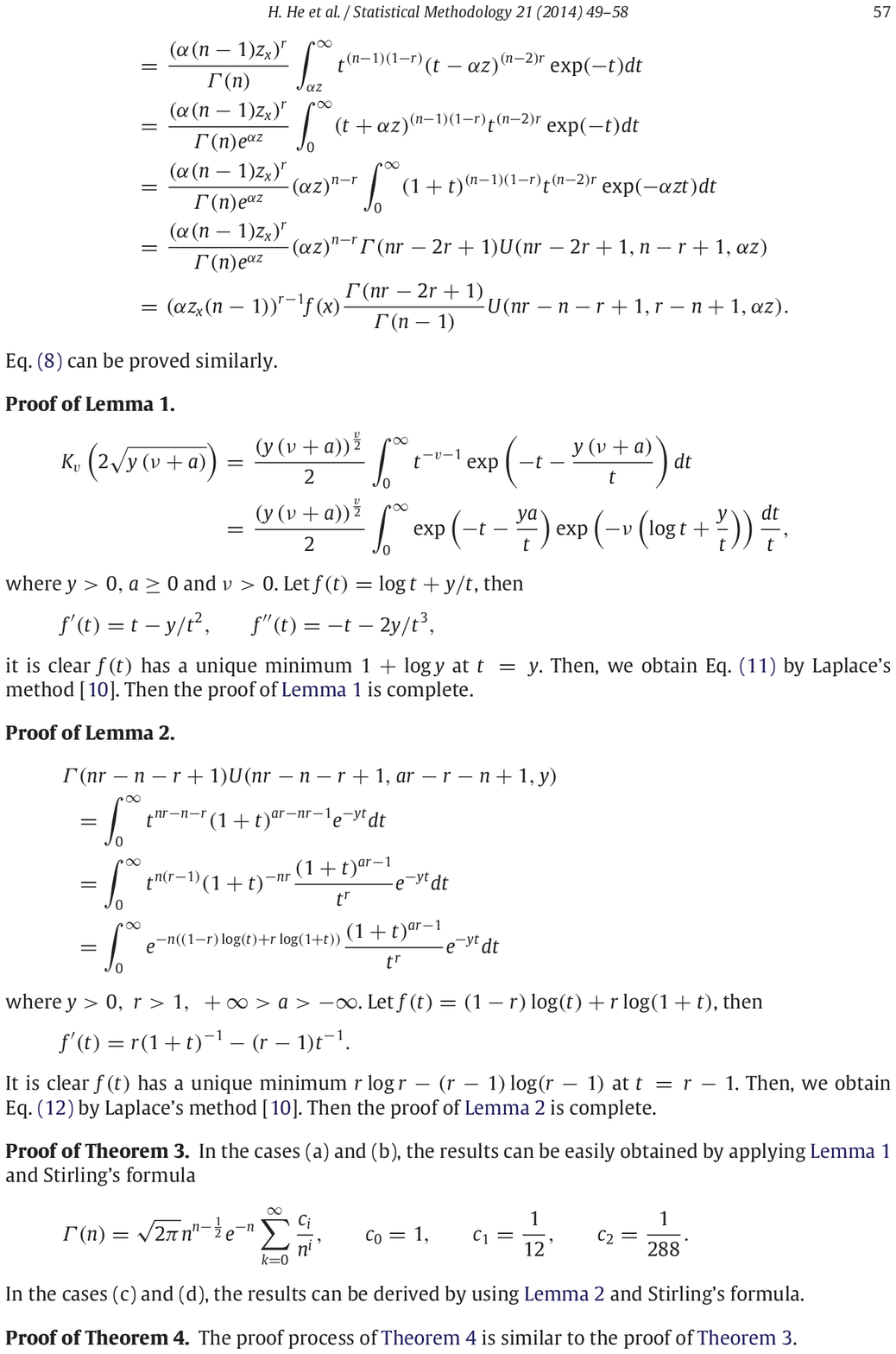}
\end{tabular}
\end{figure}
\input{epsf}
\begin{figure}
\begin{tabular}{rr}
\epsfxsize=6.5in\epsfysize=10in \epsffile{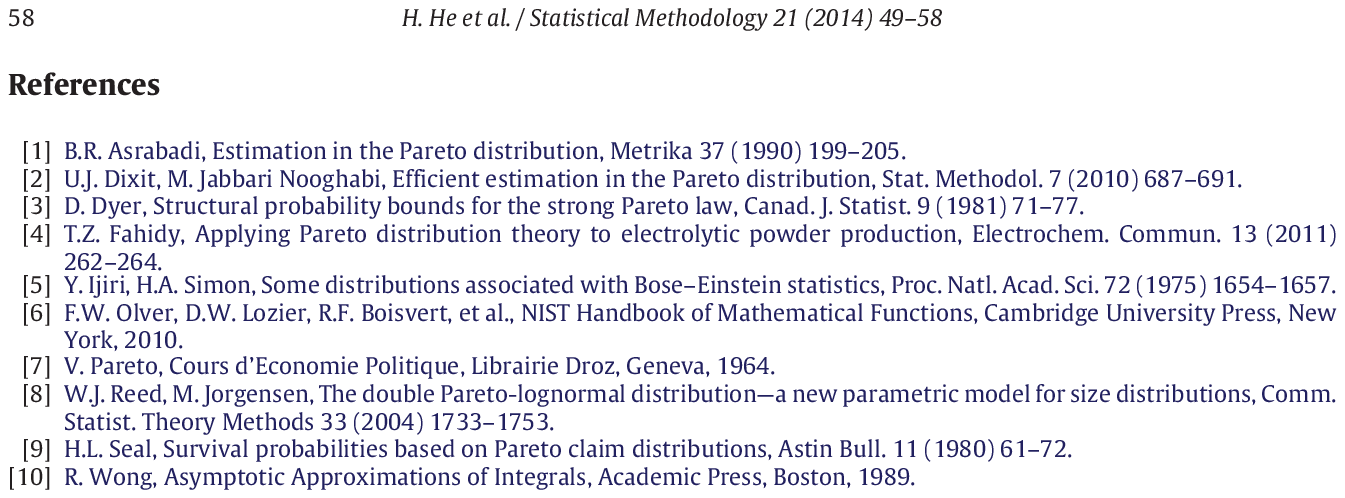}
\end{tabular}
\end{figure}

\end{document}